\documentclass[11pt]{article}

\usepackage[a4paper,margin=2cm]{geometry}
\usepackage[T1]{fontenc}
\usepackage{graphicx}
\usepackage{booktabs}
\usepackage{mathtools, amssymb}
\usepackage[dvipsnames]{xcolor}
\usepackage{subcaption}
\usepackage{algorithm, algpseudocode} 
\usepackage{siunitx}
\usepackage{makecell}
\usepackage{natbib}

\usepackage[colorlinks,bookmarksopen,bookmarksnumbered,
            citecolor=red,urlcolor=red]{hyperref}
\usepackage{cleveref}

\usepackage{setspace}
\usepackage{lineno}
\title{Splitting horizontal and vertical polynomial order in a compatible finite element discretisation for numerical weather prediction}

\author{
Daniel Witt$^{1}$\thanks{Corresponding author: D.witt@exeter.ac.uk}
\and Jemma Shipton$^{1}$
\and Thomas M Bendall$^{2}$
}

\date{}

\begin{document}

\maketitle

\noindent
$^{1}$Department of Mathematics, University of Exeter, Exeter, United Kingdom \\
$^{2}$Dynamics Research, UK Met Office, Exeter, United Kingdom

\section{Abstract}
The accurate and efficient representation of atmospheric dynamics remains a central challenge in numerical weather prediction. A particular difficulty arises from the strong anisotropy of the atmosphere, in which horizontal and vertical motions occur on very different length scales, motivating numerical discretisations that can reflect this structure. In this study, we introduce a compatible finite element discretisation of the compressible Boussinesq and compressible Euler equations in which the horizontal and vertical polynomial orders are treated independently.

The split-order discretisation is constructed using a tensor-product framework that preserves the discrete de Rham complex and associated mimetic properties. Its wave-propagation characteristics are examined through a discrete dispersion analysis that extends previous analyses to configurations with differing horizontal and vertical polynomial orders. The results show that increasing horizontal order improves the representation of gravity waves at low and intermediate wavenumbers, while increasing vertical order can degrade dispersion accuracy near the grid scale and introduce spectral gaps.

A series of idealised numerical experiments, including gravity-wave propagation, advective transport, mountain-wave flow, and a global baroclinic-wave test, is used to assess the scheme's accuracy and convergence properties. These experiments demonstrate that increasing the polynomial order in the dominant direction of motion improves convergence, and that increasing the horizontal order yields the greatest gain in accuracy under typical atmospheric conditions. The results indicate that split-order compatible finite element discretisations provide a viable alternative for controlling accuracy and numerical behaviour in atmospheric dynamical cores.

\section{Introduction}
The UK Met Office has developed a new dynamical core, LFRic, which uses a compatible finite-element discretisation \citep{adamsLFRicMeetingChallenges2019, melvinMixedFiniteelementFinitevolume2024, melvinMixedFiniteelementFinitevolume2019} that preserves many of the desirable properties of the Met Office's previous ENDGame dynamical core \citep{woodInherentlyMassconservingSemiimplicit2014, melvinInherentlyMassconservingIterative2010}.
ENDGame is a semi-implicit semi-Lagrangian finite difference model discretised on a latitude-longitude grid, with Arakawa C-grid and Charney-Phillips staggerings \citep{arakawaComputationalDesignBasic1977, charneyNUMERICALINTEGRATIONQUASIGEOSTROPHIC1953}. This discretisation choice ensures that ENDGame has good wave dispersion properties and avoids spurious pressure modes, which \cite{staniforthHorizontalGridsGlobal2012} identify as necessary and desirable features for a dynamical core. These features ensure that the dynamical core correctly captures wave propagation and key atmospheric features, such as geostrophic balance, without resorting to numerical diffusion to suppress computational noise. Despite these desirable properties, ENDGame's reliance on a latitude-longitude grid poses a significant hurdle with the increasingly parallel architecture of modern supercomputers; higher-resolution latitude-longitude grids create a communication bottleneck due to the clustering of grid cells at the poles.
Compatible finite element methods retain the desirable wave-propagation properties and eliminate spurious pressure modes whilst also facilitating non-orthogonal quasi-uniform grids that circumvent the pole-clustering problem \citep{cotterMixedFiniteElements2012, cotterCompatibleFiniteElement2023, nataleCompatibleFiniteElement2016, cotterCompatibleFiniteElement2022}.
In finite element methods, increasing the polynomial order offers an alternative to increasing resolution to improve model accuracy. Resolved motions in a finite element scheme are determined not only by the grid resolution but also by the node distribution within an element, which is governed by the polynomial order. Unlike finite difference or finite volume schemes, which rely on larger stencils for higher-order representation (such as computing the pressure gradient), higher-order finite element methods simply increase the polynomial order of representation, without the increased cell-to-cell communication associated with larger stencils.
\cite{nairComputationalAspectsScalable2009} demonstrated that higher-order finite elements still retain strong computational scaling in massively parallel architectures.

We construct two and three-dimensional finite element spaces by taking the tensor product of finite element spaces representing the horizontal discretisation with those representing the vertical (see \cref{sec: formulation}). Many finite element discretisations are constructed this way, including the compatible variety \citep{cotterFiniteElementExterior2014, arnoldDifferentialComplexesStability2006}. 
This tensor product construction provides a framework in which the horizontal and vertical orders can be changed independently. For example, \cite{guerraHighorderStaggeredFiniteelement2016} developed a general spectral element scheme for the vertical coordinate, which used higher order vertical elements to eliminate stationary modes and maintain hydrostatic balance whilst keeping the horizontal order relatively coarse. 

We are motivated to treat the horizontal and vertical directions separately due to the anisotropy of the atmosphere, where horizontal length scales are significantly larger than vertical ones. How motion at these scales is represented in numerical models is intrinsic to their accuracy and performance. \cite{skamarock_vertical_2019} showed that mesoscale-inertia gravity waves were particularly sensitive to the vertical resolution, assuming that the horizontal resolution was well resolved. These waves are crucial to the mesoscale KE spectrum in the stratosphere, which is especially important for high-resolution global forecasts. However, for coarser models, such as those used in ensemble or climate forecasts, a coarser vertical resolution is appropriate to better filter spurious flow features without seriously deforming the KE spectra \citep{burgess_troposphere--stratosphere_2013, waite_mesoscale_2009}. 
\cite{schmidtEffectsVerticalGrid2024} found results indicating that global storm-resolving models also benefited more from an increase in horizontal resolution than vertical resolution. It is therefore desirable that a model's spatial discretisation remain flexible, allowing the best combination of vertical and horizontal orders to be chosen for the particular modelling goal. This can extend to the model construction; for instance, LFRic stores data in an unstructured manner in the horizontal direction but retains vertical structure by ensuring that vertically adjacent data points are adjacent in memory \citep{adamsLFRicMeetingChallenges2019}. Therefore, changing the polynomial order in the vertical, not the horizontal, or vice versa, may have different effects on the balance between accuracy and computational performance. 

Higher-order schemes offer greater accuracy; however, they also introduce new challenges. One such challenge is the degradation of numerical dispersion relations \citep{lerouxAnalysisNumericallyInduced2007, melvinTwodimensionalMixedFiniteelement2014, cotterCompatibleFiniteElement2023}. Although compatible finite elements provide a good representation of the dispersion relation for low wavenumbers, there is a degradation at wavenumbers approaching the grid scale \citep{melvinChoiceFunctionSpaces2018, eldredDispersionAnalysisCompatible2018}. Additionally, the higher-order node distribution within an element increases the resolved motions past the grid spacing; these resolved motions can exaggerate spurious numerical noise \citep{herringtonPhysicsDynamicsCoupling2019, hannahSeparatingPhysicsDynamics2021}. This motivates our discrete dispersion analysis to ensure that such degradations are either absent or minimal.
When designing an operational dynamical core that must run within a fixed time window, we must balance attempts to improve accuracy with the resulting computational cost. In this work, we focus solely on the numerical behaviour of the split-order scheme, leaving questions of computational performance to future work.

In this paper, we develop and implement a novel split-order finite element scheme in Gusto, a compatible finite element dynamical core toolkit. To analyse its properties, we extend the dispersion analysis of \citet{melvinChoiceFunctionSpaces2018} to cases with differing polynomial orders between the vertical and horizontal directions and complement this with a series of numerical experiments assessing the accuracy of the resulting schemes. We show that the split-order discretisation improves dispersion behaviour at slow and intermediate wavenumbers but can degrade dispersion for wavenumbers approaching the grid scale. Increasing the polynomial order also introduces spectral gaps in the dispersion relation; however, we demonstrate that these do not affect numerical results in standard test cases.
The paper is organised as follows: section~\ref{sec: formulation} 
describes the construction of the split-order compatible finite element scheme. Section~\ref{sec: dispersion} presents a theoretical analysis of its dispersion properties. Section~\ref{sec: model description} outlines the numerical model and experimental setup. Results from test cases are shown in Section~\ref{sec: test cases}, followed by concluding remarks.

\section{Governing equations} \label{sec: equations}
In this section, we introduce the equation sets we use in the rest of the paper. We use the compressible Boussinesq equations for a discrete dispersion analysis since this equation set has a linear pressure gradient term following \citep{melvinChoiceFunctionSpaces2018, melvinTwodimensionalMixedFiniteelement2014}. We use the compressible Euler equations to demonstrate the accuracy and convergence of the split-order method using standard
tests commonly used in the development of weather and climate models.

\subsection{Compressible Boussinesq equations} 
The compressible Boussinesq equations in a non-rotating domain are 
\begin{subequations} \label{eqn: boussinesq group}
\begin{align}
    \frac{\partial\mathbf{u}}{\partial t} + (\mathbf{u}\cdot\nabla)\mathbf{u} + \nabla p + b\hat{\mathbf{k}} &= 0, \label{eqn:boussinesq u eqn} \\
    \frac{\partial p }{\partial t} + (\mathbf{u}\cdot\nabla)p + c_s^2 \nabla\cdot\mathbf{u} &= 0, \label{eqn:boussinesq p eqn}\\
    \frac{\partial b }{\partial t} + (\mathbf{u}\cdot\nabla)b &= 0, \label{eqn:boussinesq b eqn}
\end{align} 
\end{subequations}
where $\mathbf{u}$ is the velocity vector, $p$ is the pressure, $b$ is the 
buoyancy, $\hat{\mathbf{k}}$ is a unit vector in the vertical direction, 
and $c_s=\qty{340}{\metre\per\second}$ is the speed of sound. 

\subsection{Compressible Euler equations}
The compressible Euler equations are
\begin{subequations} \label{eqn: euler group}
\begin{align}
    \frac{\partial \mathbf{u}}{\partial t} + (\mathbf{u}\cdot \nabla)\mathbf{u} + 
    2\boldsymbol{\Omega}\times \mathbf{u} + c_p\theta\nabla\Pi +  g \hat{\mathbf{k}} &=0, \\ \label{eqn:euler u}
    \frac{\partial \rho}{\partial t} + \nabla\cdot(\mathbf{u}\rho) &= 0, \\
    \frac{\partial\theta}{\partial t} + (\mathbf{u}\cdot\nabla)\theta &= 0, \\
    \Pi^{(1-\kappa)/\kappa} &= \frac{R}{p_0}\rho\theta,  
\end{align}
\end{subequations}
where $\mathbf{u}$ is the velocity vector, $\theta$ is the potential temperature, 
$\rho$ is the density, $\Pi$ is the Exner pressure, $c_p$ is the specific heat 
capacity at constant pressure, $R$ is the gas constant, $p_0$ is the reference 
pressure, $g$ is the acceleration due to gravity, $\kappa = R/c_p$, and 
$\boldsymbol{\Omega}$ is the planet's rotation vector.

\section{Spatial Discretisation} \label{sec: formulation}
As the numerical analysis in Section \ref{sec: dispersion} uses the compressible Boussinesq equations, we present the discretisation of the linear compressible Boussinesq equations here. The discretisation of the compressible Euler equations follows the same method. Full details can be found in \citep{cotterCompatibleFiniteElement2022, bendallImprovingAccuracyDiscretisations2023}.
To discretise equation set \eqref{eqn: boussinesq group} with compatible finite elements, we define a set of spaces in which to represent the prognostic variables $(\mathbf{u}, p, b)$ and their corresponding test functions $(\boldsymbol{\psi}, \phi, \gamma)$. The spaces are chosen so that, together with the appropriate differential operators, they form a discrete de Rham complex; this ensures that our discretisation preserves the identities of continuous vector calculus. Preserving these identities ensures that spurious computational modes do not arise in the numerical solution and large-scale geostrophic flows are well represented \citep{nataleCompatibleFiniteElement2016, cotterMixedFiniteElements2012, cotterCompatibleFiniteElement2022}. 

We construct two and three-dimensional finite element spaces by taking the tensor product between a horizontal finite element and a vertical finite element. The geometry of the domain determines the form of the horizontal space, an interval in two dimensions, or a quadrilateral or triangle in three dimensions. Many finite element spaces satisfy the de Rham conditions; however, here we use the sequence $(\mathbb{V}_1, \mathbb{V}_2)$ on quadrilateral cells, where $\mathbb{V}_1$ is an order $r$ quadrilateral Raviart-Thomas element, which ensures continuity in the normal direction across the cell boundaries, and $\mathbb{V}_2$ is the required order $r-1$ discontinuous space. The vertical element is always defined on an interval.

\begin{figure}
    \centering  
    \includegraphics[scale=1]{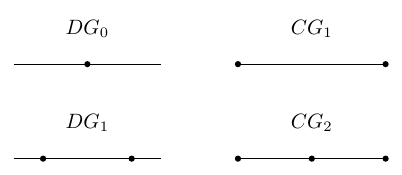}
    \caption{Diagrams of the one-dimensional finite element spaces used
    to build the two-dimensional spaces. Here, DG1 refers to discontinuous order
    1 (linear) polynomials over the cell, CG2 continuous order 2 (quadratic) polynomials,
    etc. Dots represent degrees of freedom (DoF).}
    \label{fig:interval_Diagrams}
\end{figure}

Figure \ref{fig:interval_Diagrams} illustrates the continuous ($CG^0_r$) and discontinuous ($DG^0_r$) interval finite elements used for two-dimensional domains. Here, the superscript $0$ indicates construction on an interval and $r$ denotes the polynomial order of the element. 

In addition to the spaces that form the de Rham complex, we define another space, $\mathbb{V}_{b}$, with nodes co-located with the vertical velocity nodes, giving continuity in the vertical direction. We define the buoyancy variable in $\mathbb{V}_{b}$ to mimic the Charney-Phillips grid staggering \citep{charneyNUMERICALINTEGRATIONQUASIGEOSTROPHIC1953, melvinChoiceFunctionSpaces2018, melvinMixedFiniteelementFinitevolume2024, bendallSolutionTrilemmaMoist2023}.
%Each prognostic variable and its test function is assigned to the corresponding space $\mathbf{u}, \boldsymbol{\psi} \in \mathbb{V}_1$, $p, \phi \in \mathbb{V}_2$, and $b, \gamma \in \mathbb{V}_b$, where $\boldsymbol{\psi} = (\chi, \nu)$.
In this work, we consider a set of de Rham complexes with different horizontal and vertical order. These cases are denoted throughout by the tuple ($h$, $v$), where $h$ and $v$ denote the horizontal and vertical orders in the $\mathbb{V}_2$ space, respectively. Our focus is on combinations of zeroth and first-order complexes, namely the cases $(0,0)$, $(1,0)$, $(0,1)$, and $(1,1)$, as illustrated in \cref{tab:AllSpaces}. To construct these finite element spaces in a vertical slice, we take the following tensor products,
\begin{subequations}
\begin{align}
    \mathbb{V}_{2}^{h,v} &= DG^{0, x}_{h} \otimes DG^{0, z}_{v}, \label{space:CG} \\
    \mathbb{V}_{1}^{h,v} &= CG^{0, x}_{h+1} \otimes DG^{0, z}_{v} + DG^{0, x}_{h}\otimes CG^{0, z}_{v+1}, \label{space:RT} \\
    \mathbb{V}_{b}^{h,v} &= DG^{0, x}_{h} \otimes CG^{0, z}_{v+1}.
    \label{space:VCG} 
\end{align}
\end{subequations}
\begin{table}
    \centering
    \resizebox{0.95\textwidth}{!}{
    \begin{tabular}{ccccc}
      \toprule
      Case & \multicolumn{2}{c}{$\mathbb{V}_1$} & $\mathbb{V}_2$ & $\mathbb{V}_b$ \\
      \midrule
      % First Row 
        $(0,0)$  &
        \multicolumn{2}{c}{$CG_1 \otimes DG_0 \;+\; DG_0 \otimes CG_1$} &
        $DG_0 \otimes DG_0$ &
        $DG_0 \otimes CG_1$
        \\
        &
        \hspace*{0.8em}\includegraphics{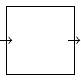} &
        \hspace*{1.4em}\includegraphics{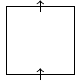} &
        \includegraphics{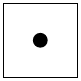} &
        \includegraphics{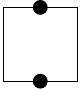}
        \\
      \midrule
      % Second Row
        $(0,1)$ &
        \multicolumn{2}{c}{$CG_1 \otimes DG_1 \;+\; DG_0 \otimes CG_2$} &
        $DG_0 \otimes DG_1$ &
        $DG_0 \otimes CG_2$
        \\
        &
        \hspace*{0.8em}\includegraphics{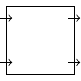} &
        \hspace*{1.4em}\includegraphics{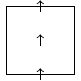} &
        \includegraphics{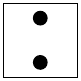} &
        \includegraphics{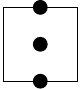}
        \\
      \midrule
      % Third Row
        $(1,0)$ &
        \multicolumn{2}{c}{$CG_2 \otimes DG_0 \;+\; DG_1 \otimes CG_1$} &
        $DG_1 \otimes DG_0$ &
        $DG_1 \otimes CG_1$
        \\
        &
        \hspace*{0.8em}\includegraphics{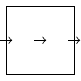} &
        \hspace*{1.4em}\includegraphics{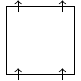} &
        \includegraphics{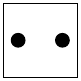} &
        \includegraphics{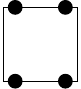}
        \\
      \midrule
      % Fourth Row
        $(1,1)$ &
        \multicolumn{2}{c}{$CG_2 \otimes DG_1 \;+\; DG_1 \otimes CG_2$} &
        $DG_1 \otimes DG_1$ &
        $DG_1 \otimes CG_2$
        \\
        &
        \hspace*{0.8em}\includegraphics{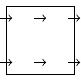} &
        \hspace*{1.4em}\includegraphics{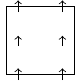} &
        \includegraphics{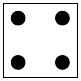} &
        \includegraphics{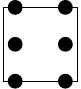}
        \\
      \bottomrule
    \end{tabular}}
    \caption{Function spaces constructed for each case. A filled circle on a cell edge indicates continuity across cell boundaries, while an arrow indicates the continuity of normal components of a vector across a boundary (but not continuity for the tangential components} 
    \label{tab:AllSpaces}
\end{table}

Linearising equations (\ref{eqn:boussinesq u eqn} -\ref{eqn:boussinesq b eqn}) around
a state of rest and multiplying each governing equation by a test function $(\boldsymbol{\psi}, \phi, \gamma) 
\in \mathbb{V}_1 \times \mathbb{V}_2 \times \mathbb{V}_b$ then integrating over 
the domain $\Omega$ yields the spatial discretisation
\begin{subequations}
\begin{align}
    \int_{\Omega}\frac{\partial u}{\partial t} \chi \;d\Omega - \int_{\Omega}p\frac{\partial \chi}{\partial x}\; d\Omega  &= 0, \; \quad \forall \chi \in \mathbb{V}_x, \label{eqn: u horizontal}\\
    \int_{\Omega}\frac{\partial w}{\partial t} \nu  \;d\Omega - \int_{\Omega}p\frac{\partial \nu }{\partial z}\;d\Omega - \int_{\Omega} b\nu \; d\Omega &= 0, \; \quad \forall \nu \in \mathbb{V}_z,  \\
    \int_{\Omega}\frac{\partial p}{\partial t} \psi \;d\Omega + c_s^2\int_{\Omega}\left( \frac{\partial u}{\partial x}+\frac{\partial w}{\partial z}\right)\psi \; d\Omega &= 0, \quad \forall \psi \in \mathbb{V}_2,\\
    \int_\Omega \frac{\partial b}{\partial t} \gamma \;d\Omega + \int_{\Omega} N^2 w\gamma \; d\Omega &= 0, \quad \forall \gamma \in \mathbb{V}_b, \label{eqn: bouyancy}
\end{align}
\end{subequations}
where $\boldsymbol{\psi} = (\chi, \nu)$. Here, we have separated the momentum equation into vertical and horizontal componentsand integrated the pressure gradient term by parts with the boundary terms vanishing due to the continuity of $\boldsymbol{\psi}\cdot\mathbf{n}$.
$\mathbb{V}^{h,v}_x$ and $\mathbb{V}^{h,v}_z$ are the components of equation \eqref{space:RT} such that $\mathbb{V}^{h,v}_1 = \mathbb{V}^{h,v}_x + \mathbb{V}^{h,v}_z$ with
\begin{align}
    \chi \in\mathbb{V}^{h, v}_x &= CG^{0, x}_{h+1} \otimes DG^{0, z}_{v} \label{space:RTX}, \\
    \nu \in \mathbb{V}^{h, v}_z &= DG^{0, x}_{h} \otimes CG^{0, z}_{v+1} \label{space:RTZ}, 
\end{align}

We now calculate a discrete dispersion relation to validate the wave propagation properties of each discretisation.

\section{Discrete dispersion analysis} \label{sec: dispersion}
A dispersion relation relates the frequency and wavenumber of waves. 
In the discrete setting, the numerical scheme determines this relation. Accurate discrete dispersion relations ensure that wave phenomena, such as inertial, gravity, and Rossby waves, are well-represented. Poor representation of such waves may result in incorrect wave speeds or computational modes. We will focus on vertical-slice test cases in which rotational forces are absent. Hence, we will not consider Rossby waves.

\subsection{Analytic dispersion relation}
\begin{figure}
    \centering
    \includegraphics[width=0.9\linewidth]{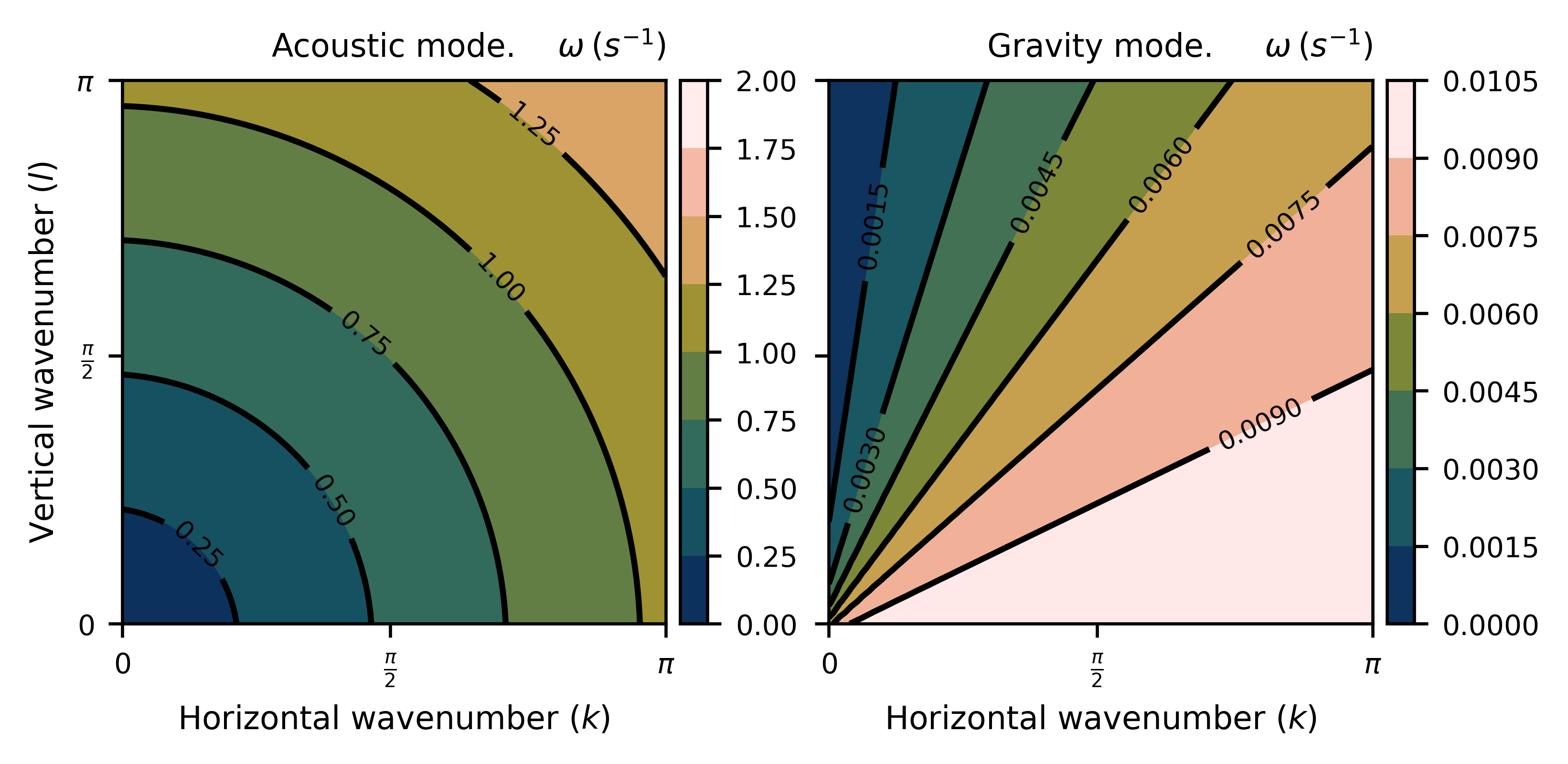}
    \caption{Analytical dispersion relation for the acoustic (left) and gravity (right)
    modes of the linear compressible Boussinesq equations.}
    \label{fig:analytical dispersion relation}
\end{figure}

We recall the analytical dispersion relation of the linearised compressible Boussinesq equations as a reference for the discrete analysis carried out in the later section. Seeking wavelike solutions of the form 
\begin{equation}
    \label{wavelike sols}
    (u, w, p, b) = (U, W, P, B)\exp{[i(kx +lz - \omega t)}].
\end{equation}
The linear compressible Boussinesq equations give the matrix vector system, 
\begin{equation}
\left( \begin{array}{cccc}
-i\omega & 0 & -ik & 0\\
0 & -i\omega & -il & -1 \\
ikc_s^2 & ilc_s^2  & -i\omega & 0 \\
0 & N^2 & 0 & -i\omega
\end{array} \right)
\left( \begin{array}{c}
U \\
W \\
P \\
B 
\end{array} \right) = 0.
\end{equation}
The system admits a unique solution only when the determinant of the matrix vanishes. Hence, 
\begin{equation}
    \omega^4 - \omega^2[(k^2 + l^2)c_s^2 + N^2] + k^2N^2c_s^2 = 0,
\end{equation}
which is a biquadratic equation for $\omega$, with solutions
\begin{equation}
    \label{eqn: modes}
    \omega^2 = \frac{1}{2}[N^2 + (k^2 + l^2)c_s^2]\left\{ 1 \pm \sqrt{1 - \frac{4k^2N^2c_s^2}{N^2 + (k^2 + l^2)c_s^2}}\right\}.
\end{equation}
The solution with a positive sign $\omega_{+}$ corresponds to acoustic waves, whilst the solution with a negative sign $\omega_{-}$ corresponds to gravity waves \citep{durranNumericalMethodsFluid2010}. Assuming $N^2 \ll k^2 c_s^2$, which is valid for most atmospheric applications, we simplify the equation to yield the two modes
\begin{equation}
    \omega_+^2 \approx N^2 + (k^2 + l^2)c_s^2, \qquad 
    \omega_-^2 \approx \frac{k^2 N^2}{\dfrac{N^2}{c_s^2} + k^2 + l^2}.
\end{equation}
These are show in figure \ref{fig:analytical dispersion relation}. 
The frequency $\omega$ is purely real, so all modes are neutrally stable and neither grow nor decay over time. Moreover, the dispersion relations are monotonic and non-zero for all $(k,l) > 0$. We now investigate the discrete dispersion relation and demonstrate that these features persist.
\subsection{Discrete dispersion relation}
Recall the discretisation of the linearised compressible Boussinesq equations(\ref{eqn: u horizontal} - \ref{eqn: bouyancy}).To construct the function spaces defined in table \ref{tab:AllSpaces}, we define a set of one-dimensional quadratic $E$, linear $F$, and constant $G$ basis functions 
across an interval $s\in[s^m, s^{m+1}]$,
\begin{subequations}
\begin{align}
    G^{m+1/2}(s) &= 1, \\
    F^m(s) &= \frac{s^{m+1} - s}{s^{m+1} - s^{m}}, \\
    F^{m+1}(s) &= \frac{s - s^{m}}{s^{m+1} - s^{m}}, \\
    E^m(s) &= \left(\frac{s-s^{m+1/2}}{s^{m}-s^{m+1/2}}\right)\left(\frac{s-s^{m+1}}{s^{m} - s^{m+1}}\right), \\
    E^{m+1/2}(s) &= \left(\frac{s-s^{m}}{s^{m+1/2}-s^{m}}\right)\left(\frac{s-s^{m+1}}{s^{m+1/2} - s^{m+1}}\right), \\
    E^{m+1}(s) &= \left(\frac{s-s^{m}}{s^{m+1}-s^{m}}\right)\left(\frac{s-s^{m+1/2}}{s^{m+1} - s^{m+1/2}}\right).
\end{align}
\end{subequations} 
Taking tensor products of these one-dimensional basis functions defines the two-dimensional basis functions. We now define the derivation of the dispersion relation for the $(1,0)$ case.The other cases follow similarly. In the $(1, 0)$ case, the basis functions are linear and quadratic for the horizontal direction and linear and constant for the vertical direction
\begin{subequations} 
\begin{align}
\chi_{1,2,3} &= \{E^{m}(x)G^{n+1/2}(z)\:,E^{m+1/2}(x)G^{n+1/2}(z)\:, E^{m+1}(x)G^{n+1/2}(z)\:\}, \label{eqn: horizontal basis} \\  
(\nu, \gamma)_{1,2,3,4} &= \{F^{m}(x)F^{n}(z)\:,F^{m+1}(x)F^{n}(z)\:, F^{m}(x)F^{n+1}(z)\:, F^{m+1}(x)F^{n+1}(z)\: \}, \\
\psi_{1, 2} &= \{F^{m}(x)G^{n+1/2}(z)\:, F^{m+1}(x)G^{n+1/2}(z)\}. \label{eqn: DG basis}
\end{align}
\end{subequations} 
The prognostic variables restricted to a single element $[x^m, x^{m+1}]\times[z^m, z^{m+1}]$ are
\begin{subequations}
\begin{align}
    \hat{u}(x,z) &= u_{(0,+)}^{m, n+1/2}\chi_1 + u_{(+,+)}^{m+1/2, n+1/2}\chi_2 + u_{(0,+)}^{m+1, n+1/2}\chi_3, \label{eqn: u expression}\\
    \hat{w}(x, z) &= v_{(+,0)}^{m, n}\nu_1 +v_{(+,0)}^{m+1, n}\nu_2 +v_{(-,0)}^{m, n+1}\nu_3 + v_{(-,0)}^{m+1, n+1}\nu_4, \label{eqn: w expression}\\
    \hat{p}(x, z) &= p_{(+,+)}^{m, n+1/2}\psi_1 +p_{(-,+)}^{m+1, n+1/2}\psi_2  \label{eqn: p expression}\\
    \hat{b}(x, z) &= b_{(+,0)}^{m, n}\gamma_1 +b_{(+,0)}^{m+1, n}\gamma_2 +b_{(-,0)}^{m, n+1}\gamma_3 + b_{(-,0)}^{m+1, n+1}\gamma_4, \label{eqn: b expression}
\end{align}
\end{subequations}
where the subscripts $(0, -, +)$ denote whether the field is continuous $(0)$ or discontinuous $(+,-)$ at this DoF, the $(+,-)$ denotes which side of the facet the discontinuous DoF belongs. It is helpful to number the DoFs for each cell such that,
\begin{subequations}
    \begin{align}
        (u_1, u_2, u_3) &= (u^{m,n+1/2},u^{m+1/2,n+1/2},u^{m+1,n+1/2}), \\
        (w_1, w_2, w_3, w_4) &= (w^{m,n}, w^{m+1,n}, w^{m,n+1}, w^{m+1,n+1}), \\
        (p_1, p_2) &= (p^{m,n+1/2}, p^{m+1,n+1/2}),  \\
        (b_1, b_2, b_3, b_4) &= (b^{m,n}, b^{m+1,n}, b^{m,n+1}, b^{m+1,n+1}).
    \end{align}
\label{eqn:number_guide}
\end{subequations}
Figure \ref{fig:DoF locations 2} shows the locations of these nodes (node locations of the (0, 1) can be found in appendix \ref{A:(0,1)locs}).
\begin{figure}
    \centering
    \includegraphics[scale=1]{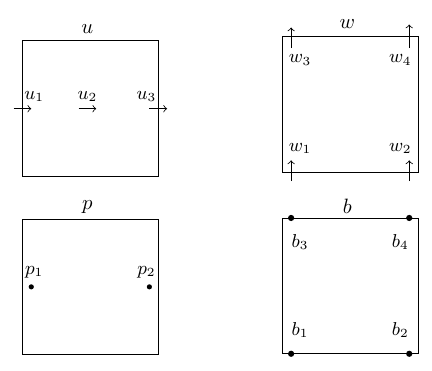}
    \caption{Locations of degrees of freedom for the $(1,0)$ case for the 
             horizontal velocity field $(u)$ shown in the top left, followed by the 
             vertical velocity field $(w)$ in the top right, the pressure field $(p)$ in the bottom left, and the buoyancy field $(b)$ in the bottom right. The node numbering can be related to node locations using equations \eqref{eqn:number_guide}. Nodes on a cell edge indicate continuity between cells, and an arrow across an edge indicates normal continuity.}
    \label{fig:DoF locations 2}
\end{figure}
Substituting the variable expansions (\ref{eqn: u expression} - \ref{eqn: b expression}) into the discretised equations (\ref{eqn: u horizontal} - \ref{eqn: bouyancy}), gives an equation for each distinct DoF. For example, consider solving equation \eqref{eqn: u horizontal} on a generic cell $(x,z)\in [x^m, x^{m+1}]\times[z^n, z^{n+1}]$, this yields the system of equations
\begin{subequations} 
\label{eqn: 3 equation set}
\begin{align}
    &\int^{x^{m+1}}_{x^m}\int^{z^{n+1}}_{z^n}\frac{\partial\hat{u}}{\partial t}\chi_1 \;dxdz + \int^{x^{m+1}}_{x^m}\int^{z^{n+1}}_{z^n}\hat{p}\frac{\partial \chi_1}{\partial x} \;dxdz= 0, \\ 
    &\int^{x^{m+1}}_{x^m}\int^{z^{n+1}}_{z^n}\frac{\partial\hat{u}}{\partial t}\chi_2 \;dxdz + \int^{x^{m+1}}_{x^m}\int^{z^{n+1}}_{z^n}\hat{p}\frac{\partial \chi_2}{\partial x} \;dxdz= 0, \\ 
    &\int^{x^{m+1}}_{x^m}\int^{z^{n+1}}_{z^n}\frac{\partial\hat{u}}{\partial t}\chi_3 \;dxdz + \int^{x^{m+1}}_{x^m}\int^{z^{n+1}}_{z^n}\hat{p}\frac{\partial \chi_3}{\partial x} \;dxdz= 0.    
\end{align}
\end{subequations}
As in the analytical case, we assume that all variables are wavelike,
\begin{equation}
    (u, w, p, b)^{m, n} \equiv (U, W, P, B)\exp[i(kx^m + lz^n - \omega t)], \label{eqn:u_fourier}
\end{equation}
where the $(m,n)$ superscripts define the location of the DoF as in \eqref{eqn:number_guide}. Exploiting the wave structure, we reduce the three-equation system \eqref{eqn: 3 equation set} to two equations. The continuity of the field at the DoF allows us to relate the field values at the left and right edge DoFs, leaving only two unique DoFs per cell. Explicitly, this is done by first relating the DoF $u_3$ to $u_1$ via, 
\begin{equation}
    u_3 = u_1e^{ik\Delta x},
\end{equation}
such that,
\begin{equation}
    \hat{u}(x,z) = u_1\chi_1 + u_2\chi_2 + u_1e^{ik\Delta x}\chi_3. \label{eqn: u expression related}
\end{equation}
Having related $u_3$ to $u_1$, we next relate their corresponding equations. 
Each $m$ is replaced by $m-1$ in equation \eqref{eqn: u expression}, 
this results in a `shifting' of the equation by a factor of $e^{-ik\Delta x}$ to 
yield,
\begin{equation}
    \hat{u}_{right} = u_1\chi_1e^{(-ik\Delta x)} + u_2\chi_2e^{(-ik\Delta x)} + u_1\chi_3.
\end{equation}
We then add this `shifted' $u_3$ equation to the equation for $u_1$ to get an 
expansion centred on $u_1$,
\begin{align}
    \hat{u}_{center} &= \hat{u} + \hat{u}_{right} \nonumber \\ 
    &=  u_1\chi_1e^{(-ik\Delta x)} + u_2\chi_2e^{(-ik\Delta x)} + u_1(\chi_3 + \chi_1) + u_2\chi_2 + u_1\chi_3e^{ik\Delta x}. \label{eqn:u_center_equations}
\end{align}
We repeat this process for the pressure expansion equation (\ref{eqn: p expression}). 
Unlike the horizontal velocity expansion, there is no continuity to relate nodes 
to each other, hence
\begin{align}
    \hat{p}_{center} &= \hat{p} + \hat{p}_{right} \nonumber \\
                     &= p_1\psi_1e^{(-ik\Delta x)} + p_2\psi_2e^{(-ik\Delta x)} + p_1\psi_1 + p_2\psi_2. \label{eqn:p_center_equations}
\end{align}
The system of equations in \eqref{eqn: 3 equation set} has been reduced to a set of two equations,
\begin{subequations}
\label{eqn: 2 equation set}
\begin{align}
    &\int^{x^{m+1}}_{x^m}\int^{z^{n+1}}_{z^n}\frac{\partial\hat{u}_{center}}{\partial t}\chi_1 \;dxdz + \int^{x^{m+1}}_{x^m}\int^{z^{n+1}}_{z^n}\hat{p}_{center}\frac{\partial \chi_1}{\partial x} \;dxdz= 0, \\ 
    &\int^{x^{m+1}}_{x^m}\int^{z^{n+1}}_{z^n}\frac{\partial\hat{u}}{\partial t}\chi_2 \;dxdz + \int^{x^{m+1}}_{x^m}\int^{z^{n+1}}_{z^n}\hat{p}\frac{\partial \chi_2}{\partial x} \;dxdz= 0.     
\end{align}
\end{subequations}
The same logic can be applied to the vertical velocity and buoyancy equations, as they are continuous
in the vertical direction. We follow the same procedure of first relating $(w, b)_3$
to $(w,b)_1$ and $(w, b)_4$ to $(w,b)_2$ via,
\begin{align*}
    (w_3, b_3) = (w_1, b_1)e^{il\Delta z}, \\
    (w_4, b_4) = (w_2, b_2)e^{il\Delta z}.
\end{align*}
Next, we replace each $n$ with $n-1$ in equations (\ref{eqn: w expression}) and (\ref{eqn: b expression})
for the equations at nodes $(w, b)_{3,4}$ to `shift' each equation by a factor of $e^{il\Delta z}$,
and then add them to the corresponding equations for $(w,b)_{1,2}$.
This procedure yields a system of eight equations, two for each
variable. Writing this in matrix form yields the system
\begin{equation}
    \label{eqn: block matrix form}
    \left( \begin{array}{cccc}
     -i\omega\mathbf{M}_u & 0 & \mathbf{D}^u_x & 0\\
    0 &  -i\omega\mathbf{M}_w & \mathbf{D}_z^w & -\mathbf{Q} \\
    c_s^2 \mathbf{D}_x^p & c_s^2 \mathbf{D}_z^p & -i\omega\mathbf{M}_p & 0 \\
    0 & N^2\mathbf{Q}^T & 0 & -i\omega\mathbf{M}_b
    \end{array} \right)
    \left( \begin{array}{c}
    \mathbf{\bar{U}} \\
    \mathbf{\bar{W}} \\
    \mathbf{\bar{P}} \\
    \mathbf{\bar{B}} 
    \end{array} \right) = 0,
\end{equation}
where $(\mathbf{\bar{U}}, \mathbf{\bar{W}}, \mathbf{\bar{P}}, \mathbf{\bar{B}})$,
are vectors of the amplitude factors for the basis functions of each prognostic variable such as $\mathbf{\bar{U}} = (U_1, U_2)^T$ etc. 
$M_{u, w, p, b}$, $D_{x, z}^{u, w, p}$ and $Q$ are all $2\times 2$ matrices that represent the action of mass matrices, discrete derivative matrices,
and projection matrices, respectively (details of these matrices can be found in appendix \ref{A:dispersion_matricies}).
Splitting equation (\ref{eqn: block matrix form}) into its symmetric and non-symmetric components 
yields 
\begin{equation}
    \label{eqn:split matrix}
    -i\omega\mathcal{M}\hat{\mathbf{x}} + \mathcal{S}\hat{\mathbf{x}} = 0,
\end{equation}
where
\begin{equation}
    \mathcal{M} = \left(\begin{array}{cccc}
        \mathbf{M}_u & 0 & 0 & 0 \\
        0 & \mathbf{M}_w & 0 & 0 \\
        0 & 0 & \mathbf{M}_p & 0 \\
        0 & 0 & 0 & \mathbf{M}_b \\
    \end{array}
    \right),
\end{equation}
\begin{equation}
    \mathcal{S} = \left(\begin{array}{cccc}
        0 & 0 & \mathbf{D}_x^u & 0 \\
        0 & 0 & \mathbf{D}_z^w & \mathbf{Q} \\
        \mathbf{D}_x^p & \mathbf{D}_x^p & 0 & 0 \\
        0 & \mathbf{Q}^T & 0 & 0 \\
    \end{array}
    \right),
\end{equation}
and $\hat{\mathbf{x}}= [U_1, U_2, W_1, W_2, P_1, P_2, B_1, B_2]^T$. As $\mathcal{M}$
is a block diagonal matrix of mass matrices, we know that its inverse exists such that equation \eqref{eqn:split matrix}
can be rearranged to give,
\begin{equation}
    \label{eqn:eigenvalue problem}
    \omega\hat{\mathbf{x}} = -i\mathcal{M}^{-1}\mathcal{S}\hat{\mathbf{x}} = \mathcal{F}\hat{\mathbf{x}}.
\end{equation}
We have now formed an eigenvalue problem which, for a given wavenumber pair $(k,l)$,
can be solved to yield eigenvalues $\omega$ corresponding to the linear acoustic and gravity modes.
\subsubsection{Wave mode allocation}
As the matrix $\mathcal{F}$ has dimensions $8\times 8$, we extract eight 
eigenvalue solutions from equation \eqref{eqn:eigenvalue problem}. These solutions correspond
to the sets, 
\begin{equation}
    \pm\{\omega_1^{+}, \omega_2^{+}\} \qquad \pm\{\omega_1^{-}, \omega_2^{-}\},
\end{equation} 
where $\omega^{+}$ represents the acoustic waves and $\omega^{-}$ the gravity waves. Taking the positive values of these sets, we now have two solutions for each regime instead of the single solution in the analytical case. As in \citep{melvinChoiceFunctionSpaces2018, melvinTwodimensionalMixedFiniteelement2014} we argue that this extra solution corresponds to an aliased wave from a higher wavenumber, supported by the increased finite element order in the horizontal direction. Specifically, we say one solution belongs to the initial wavenumber space $(k\Delta x, l\Delta z) \in (0,\pi)\times(0,\pi)$and the other solution corresponds to an $x$-extended wavenumber space $(k\Delta x, l\Delta z) \in (\pi,2\pi)\times(0,\pi)$. We must now determine which solutions belong to which space. To do this, we split the wavenumber space into sectors corresponding to the wavenumber extension and compute all possible permutations of $\omega$ for each sector. We measure the error between the discrete and analytic dispersion relations for each permutation of $\omega$, using the $L_2$ norm, and choose the permutation with the lowest error. To calculate the dispersion relation for the $(0,1)$ and $(1,1)$ cases we follow the same procedure as the $(1,0)$ case, with the $(0, 1)$ case being extended in the $z$ direction and $(1,1 )$ in both directions.

We examine how closely our discrete dispersion relations match the analytical dispersion relation plotted in figure \ref{fig:analytical dispersion relation}. 
\begin{figure}
    \centering
    \includegraphics[width=\linewidth]{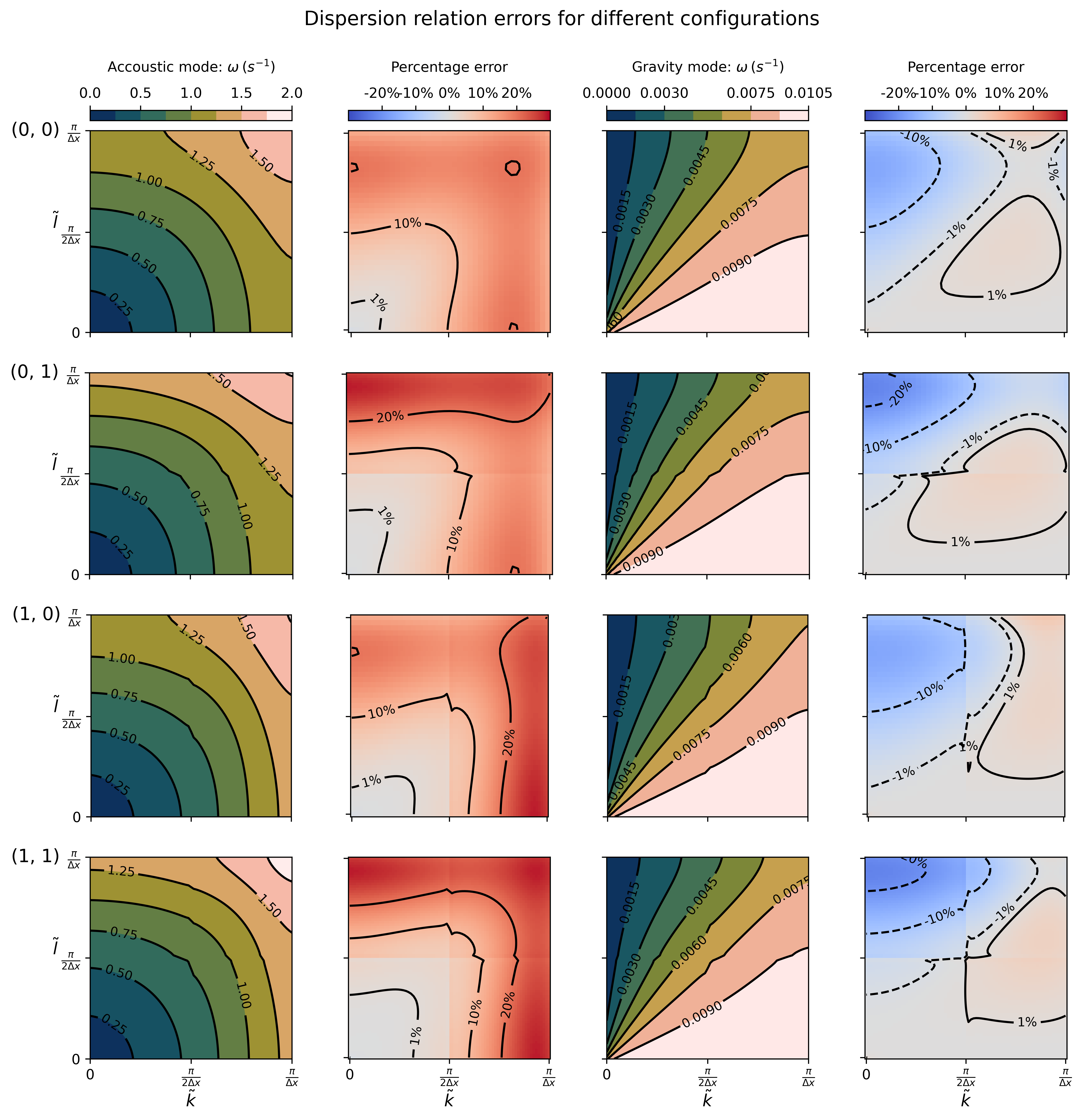}
    \caption{Discrete dispersion relations for the $(1, 1)$, $(1, 0)$ and $(0, 1)$ 
    finite element cases, for the scaled wavenumbers $\tilde{k}$, 
    $\tilde{l}$ $\in[0, \frac{2\pi}{\Delta x}] \times [0, \frac{2\pi}{\Delta z}]$, 
    where we arbitarity take $(\Delta x, \Delta z) = (1000, 1000)\unit{\meter}$. 
    To the right of each mode is a plot showing the percentage error between the discrete and 
    analytical dispersion relations.}
    \label{fig:Dispersion and error plot}
\end{figure}
Figure \ref{fig:Dispersion and error plot} shows the discrete dispersion relation for each of the four cases. To allow for comparison to the analytical dispersion relation, we use scaled wavenumbers  
\begin{equation}
    (\tilde{k},\tilde{l}) \in\left[0, \frac{2\pi}{\Delta x}\right] \times \left[0, \frac{2\pi}{\Delta z}\right],
\end{equation} 
In general, all cases provide a good representation of the analytical dispersion relation, with a few noted inaccuracies. Each case overestimates 
the acoustic mode as $(\tilde{k},\tilde{l})$ approaches $(\frac{2\pi}{\Delta x}, \frac{2\pi}{\Delta z})$, whilst the gravity mode is underestimated for waves with high $\tilde{l}$ but low $\tilde{k}$. Although we have derived the discrete dispersion relations for gravity and acoustic waves, we now focus solely on gravity waves, as acoustic waves are not important for weather and climate modelling and are typically damped or filtered by numerical schemes. When comparing different orders, we find that increasing order in either direction introduces spectral gaps in the dispersion relation along the lines $\tilde{k} = \pi / \Delta x$ and $\tilde{l} = \pi / \Delta z$, in addition higher order cases represent larger wavenumbers less accurately however improve the representation of wavenumbers  $\tilde{k}, \tilde{l} < \frac{\pi}{2\Delta x}$.
\begin{table}
\centering
\setlength{\tabcolsep}{4pt}
{\renewcommand{\arraystretch}{1.6}
\begin{tabular}{lrrr}
\hline
Case & \makecell{Normalised $L_2$\\ error norm (\qty{e-4})} & Max error & Min error \\
\hline
(0, 0) & 1.42719 & 0.0278812 & 1.89524e-10 \\
(0, 1) & 1.63875 & 0.0496205 & 1.30731e-10 \\
(1, 0) & 1.2243  & 0.0279968 & 2.92222e-11 \\
(1, 1) & 1.36044 & 0.0497644 & 3.9529e-12  \\
\hline
\end{tabular}}
\caption{Global error measured for the gravity mode. Increasing horizontal order improves the representation of the gravity mode; however, increasing the vertical order degrades it. We do, however, generally see a larger maximum error with increasing order.}
\label{tab:grav_error_statistics}
\end{table}
Table \ref{tab:grav_error_statistics} shows some error statistics for the gravity mode (acoustic mode statistics can be found in appendix \ref{A:accoustic_stats}). These statistics, combined with the results presented in figure \ref{fig:Dispersion and error plot} indicate that the vertical order is the largest contributor to overall error. The $(0, 1)$ case has the most significant maximum error and performs worse than $(0, 0)$ in the $L_2$ norm. In contrast, increasing the horizontal order, such as in $(1, 0)$ improves the $L_2$ norm error with only a marginal cost to the most extreme errors. The order $(1, 1)$ case has a worse performance when compared to the $(1, 0)$ case, but is still better than the $(0, 1)$ and $(0, 0)$ cases for the global norm, however, does possess a maximal error comparable to the $(0, 1)$ case. 
Having identified some potential shortcomings from the dispersion relation, specifically, the propagation of near grid-scale waves, we now wish to see if these issues materialise in the numerical model. 
\section{Model formulation}\label{sec: model description}
Our model is implemented in Gusto, the dynamical core toolkit built on the Firedrake finite element library. Firedrake provides an automatic code generation for the implementation of finite element methods to solve partial differential equations \citep{FiredrakeUserManual}. Gusto provides a suite of stable time stepping, transport and physics schemes for discretising and modelling both the shallow water equations and the compressible Euler equations, facilitating rapid prototyping of novel algorithms within the compatible finite element framework.

\subsection{Time discretisation}
In all test cases, the prognostic variables are evolved in time using a semi-implicit quasi-Newton time-stepping scheme. The next-generation UK Met Office model influenced this choice; a detailed description of the scheme is found in \citep{woodInherentlyMassconservingSemiimplicit2014, melvinMixedFiniteelementFinitevolume2019, bendallRecoveredSpaceAdvection2019, hartneyExploringFormsMoist2025}. Within the semi-implicit quasi-Netwon scheme, we apply an explicit transport scheme to each of the prognostic variables. In this paper, we use the SSPRK3 scheme for all variables \citep{izzo_strong_2022} and pseudo code outling the timestepping precedure can be found in appendix \ref{A:pseudo-code}.

Second-order transport is a minimum requirement for accurate weather and climate models \citep{staniforthHorizontalGridsGlobal2012}. For the $(1, 1)$ case, we employ an upwind scheme to evaluate the transport terms. However, when the model configuration includes lowest-order finite elements, we cannot achieve second-order accuracy through upwinding alone. To overcome this, we make use of the recovered space finite element transport scheme from 
\citep{bendallRecoveredSpaceAdvection2019,bendallImprovingAccuracyDiscretisations2023}.

The recovered space transport scheme achieves second-order transport by computing the transport step in a higher-order space (the embedding space) before projecting the transported field back to the original space. As long as the projection operator is second-order accurate, the entire scheme is second-order accurate. For the $(0,0)$ case \citep{bendallImprovingAccuracyDiscretisations2023} describes a sequence of operators and space that satisfy second order accuracy; however, we now must extend this to the $(1, 0)$ and $(0, 1)$ cases by defining the respective embedding spaces.

\begin{table}
    \centering
    \setlength{\tabcolsep}{2pt}
    \begin{tabular}{ccccc}
        \hline
        Variable & \multicolumn{2}{c}{Case $(1, 0)$} & \multicolumn{2}{c}{Case $(0, 1)$}  \\ \hline
                 & $V_0$      &   $V_E$               &   $V_0$    &    $V_E$ \\ \hline
        $\rho$   & 
        \includegraphics{Figures/table_plots/rho_initial_1_0.pdf} &
        \includegraphics{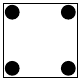} &
        \includegraphics{Figures/table_plots/rho_initial_0_1.pdf} &
        \includegraphics{Figures/table_plots/rho_recovered.pdf}  \\ \hline

        $u$ &
        \includegraphics{Figures/table_plots/u_initial_1_0.pdf} &
        \includegraphics{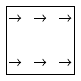}  &
        \includegraphics{Figures/table_plots/u_initial_0_1.pdf} &
        \includegraphics{Figures/table_plots/u_embed_both.pdf}  \\ \hline

        $w$ &
        \includegraphics{Figures/table_plots/w_initial_1_0.pdf} &
        \includegraphics{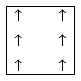}  &
        \includegraphics{Figures/table_plots/w_initial_0_1.pdf} &
        \includegraphics{Figures/table_plots/w_embed_both.pdf}  \\ \hline

        $b$ or $\theta$ &
        \includegraphics{Figures/table_plots/theta_base.pdf} &
        \includegraphics{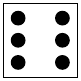} &
        \includegraphics{Figures/table_plots/theta_initial_0_1.pdf} &
        \includegraphics{Figures/table_plots/theta_embedded.pdf}  \\ \hline    
    \end{tabular}
    
    \caption{A table illustrating the initial and embedded spaces used in the $(1, 0)$
             and $(0, 1)$ cases, here $u$ and $w$ indicate the horizontal and vertical
             components of the velocity vector $\mathbf{u}$. Nodes indicate the location of DoFs within the space, with nodes placed on a cell edge indicating continuity, arrows indicate a vector value DoF and arrows across a cell edge indicate normal continuity across that edge.}
    \label{tab:recovery spaces}
\end{table}

Table \ref{tab:recovery spaces} shows the initial space $V_0$ of each variable and the embedding $V_E$ space in which the transport occurs. The field is reconstructed in the embedding space via a continuous version of the higher order space using recovery operators described in \citep{bendallImprovingAccuracyDiscretisations2023}. This method 
replaces the DG upwind method whenever a lowest-order polynomial is present.

\section{Numerical tests} \label{sec: test cases}
We now explore the accuracy and convergence properties of the split-order 
schemes over a suite of numerical test cases. 
When comparing numerical results from different order schemes, we keep the distance between DoFs constant by adjusting the cell size to account for the additional DoFs in the higher-order representation. Details of the resolutions used can be found in appendix \ref{A:resolutions}). Throughout this section, the model setup will be described in the context of the $(1, 1)$ case unless otherwise stated.

\subsection{Skamarock and Klemp gravity wave} \label{sec: gravity wave}
\begin{figure}
    \centering
    \includegraphics[scale=0.45]{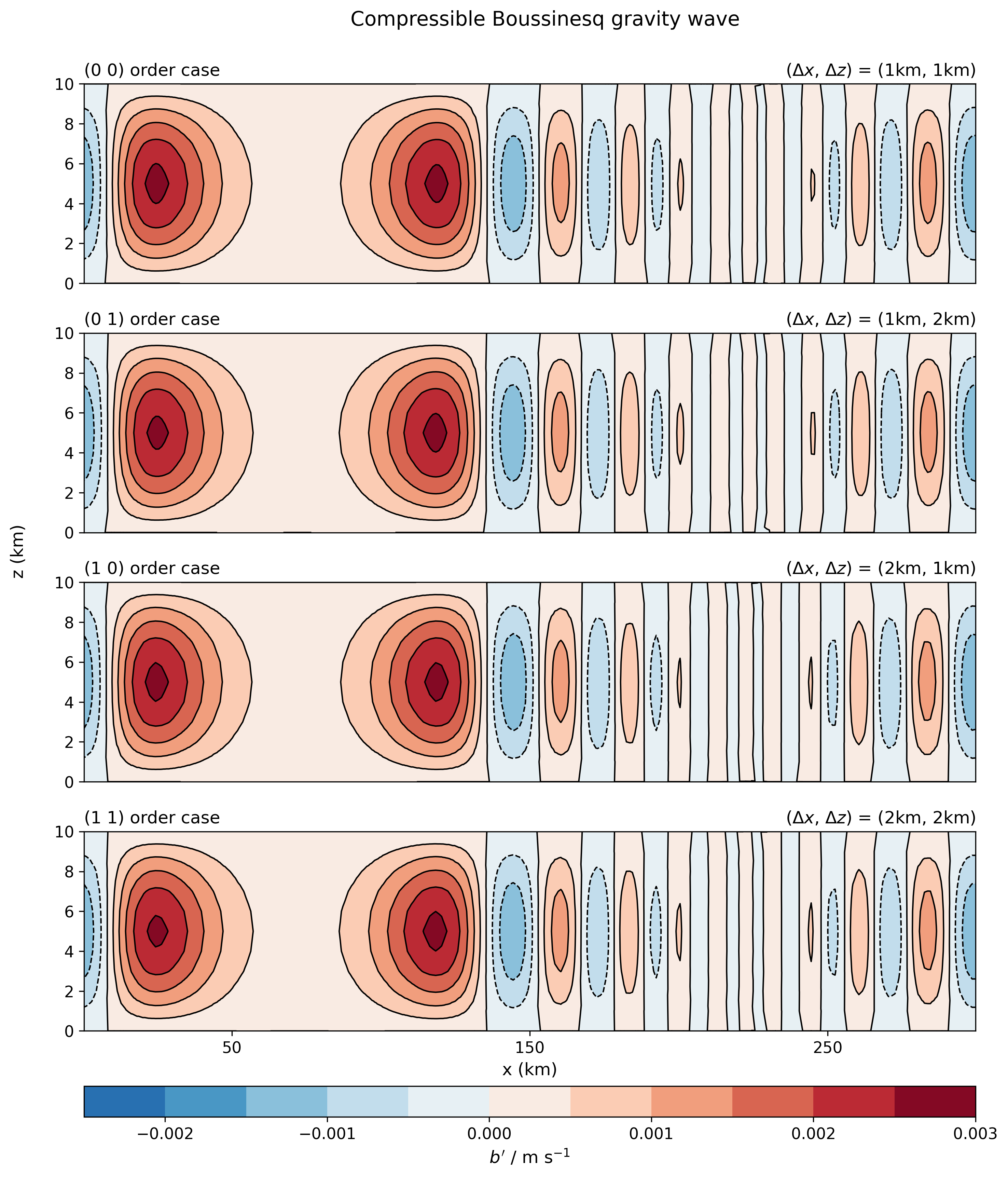}
    \caption{Buoyancy perturbation of the Skamarock and Klemp gravity wave for the 
    compressible boussinesq equations, test case at time t=3600s. Contours are between
    $[-3, 3]\times\SI{e-3}{\meter\per\square\second}$ with a spacing of 
    $0.5\times\SI{e-3}{\meter\per\square\second}$ and are the same for each case.}
    \label{fig:buoyancy plot lower}
\end{figure}

To support the theoretical findings of the discrete dispersion analysis, we consider a wave propagation benchmark of the non-hydrostatic gravity wave of \citep{skamarockEfficiencyAccuracyKlempWilhelmson1994}. This test examines the horizontal propagation of gravity waves in a two-dimensional Cartesian vertical slice domain $[0, L]\times[0, H]$ km. We initialise the test with a balanced, thermally stratified base state and introduce a buoyancy perturbation
\begin{equation}
    b' = b_0\sin{\left(\frac{\pi z}{H}\right)}\left[1 + \frac{x-x_c}{a^2}\right].
\end{equation}
Where $b_0 = \qty{0.001}{\meter\per\square\second}$, $H=\qty{10}{\kilo\meter}$, 
$a=\qty{5}{\kilo\meter}$ and centred on $x_c = \frac{L}{2}\unit{\kilo\meter}$ where $L=\qty{300}{\kilo\meter}$. 
This perturbation triggers a gravity wave that propagates along the domain. 
The integration is carried out for $t=\qty{3600}{\second}$ in time steps of $\Delta t = \qty{1.2}{\second}$.
For the order $(1,1)$ case the resolution is $(\Delta x, \Delta z)=(\qty{2000}{\meter},\qty{2000}{\meter})$ 
and the other cases have suitably scaled resolutions to ensure that each case 
has the same total number of DoFs and has the same effective resolution. 
Figure \ref{fig:buoyancy plot lower} shows the final buoyancy field, and we see 
that all cases show correct wave propagation.
We also carry out a convergence test: we compute the $L_2$ error norm by comparing simulation with $\Delta x \in [450, 500, 600, 700, 750, 800, 850, 900]\unit{\meter}$ against a high-resolution reference solution computed with $\Delta x = \qty{120}{\meter}$, which we treat as the truth. We keep vertical resolution $\Delta x = \qty{1000}{\meter}$, and time step $\Delta t = \qty{1.2}{s}$ constant for all simulations. 
\begin{figure}
    \centering
    \includegraphics[scale=0.5]{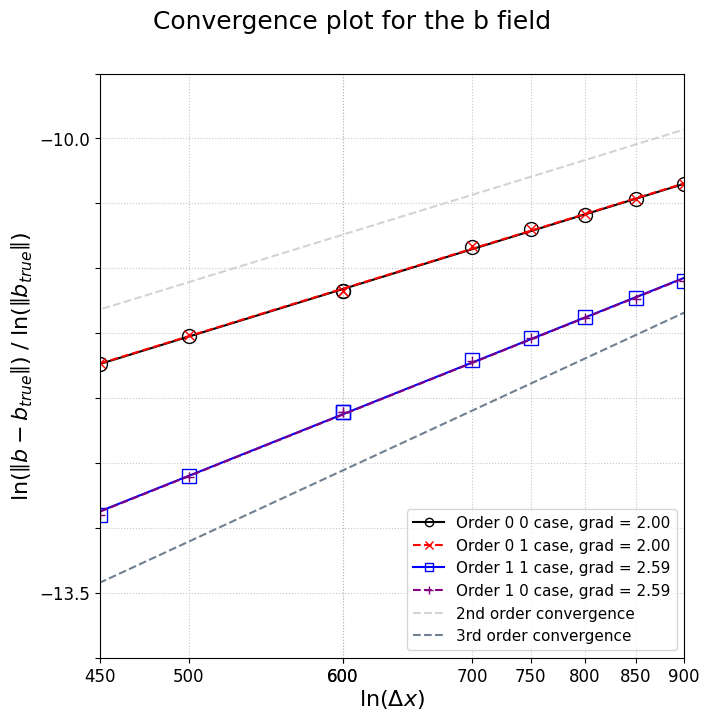}
    \caption{Convergence plot of the buoyancy field. Two comparison lines are also
    shown, to demonstrate second and third order convergence.}
    \label{fig:Gw convergence}
\end{figure}
The results presented in figure \ref{fig:Gw convergence} demonstrate an order of convergence of at least two for all cases, and the $(1,0)$ and $(1,0)$ cases have a convergence rate approaching third order, the largest order permitted by our time-stepping and transport schemes.

These results suggest that the horizontal order is dominating the accuracy for this test case. This result is expected, since the gravity wave test case consists of predominantly horizontal motion, extensive testing has found that vertical resolution generally plays little role in the convergence or accuracy of the solution.

\subsection{Travelling Vortex} \label{sec: travelling vortex}
The travelling vortex test case first proposed in \citep{kadiogluFourthorderAuxiliaryVariable2008}, formulated in line with \citep{chewOnestepBlendedSoundproofcompressible2021, kowalczykCompatibleFiniteElements2025},
consists of a smooth vortex in a doubly periodic two-dimensional domain of dimensions $[0,10]\times[0,10]\:\unit{\kilo\meter}$ with no Coriolis or gravitational forces. We initialise the vortex with a smooth radially symmetric density bump placed at the domain's centre $(5, \: 5)\unit{\kilo\meter}$.

We are interested in how each order case represents scenarios in which horizontal or vertical motion dominate. Accordingly, we prescribe three advecting velocities: a purely horizontal flow
$\mathbf{u} = (100, 0)\unit{\meter\per\second}$, a purely vertical flow
$\mathbf{u} = (0, 100)\unit{\meter\per\second}$, and a diagonal flow
$\mathbf{u} = (100, 100)\unit{\meter\per\second}$.

We run each simulation to $t = \qty{100}{\second}$ with a time step of
$\Delta t = \qty{0.1}{\second}$, ensuring that the vortex undergoes one complete traversal of the domain. The spatial resolution of the (1, 1) case is $(200, 200)\unit{\meter}$ and other cases are scaled appropriately.
For each case, we compute the $L_2$ error between the initial and final density fields and plot it as a function of the overall cell size.
\begin{figure}
    \centering
    \includegraphics[width=\linewidth]{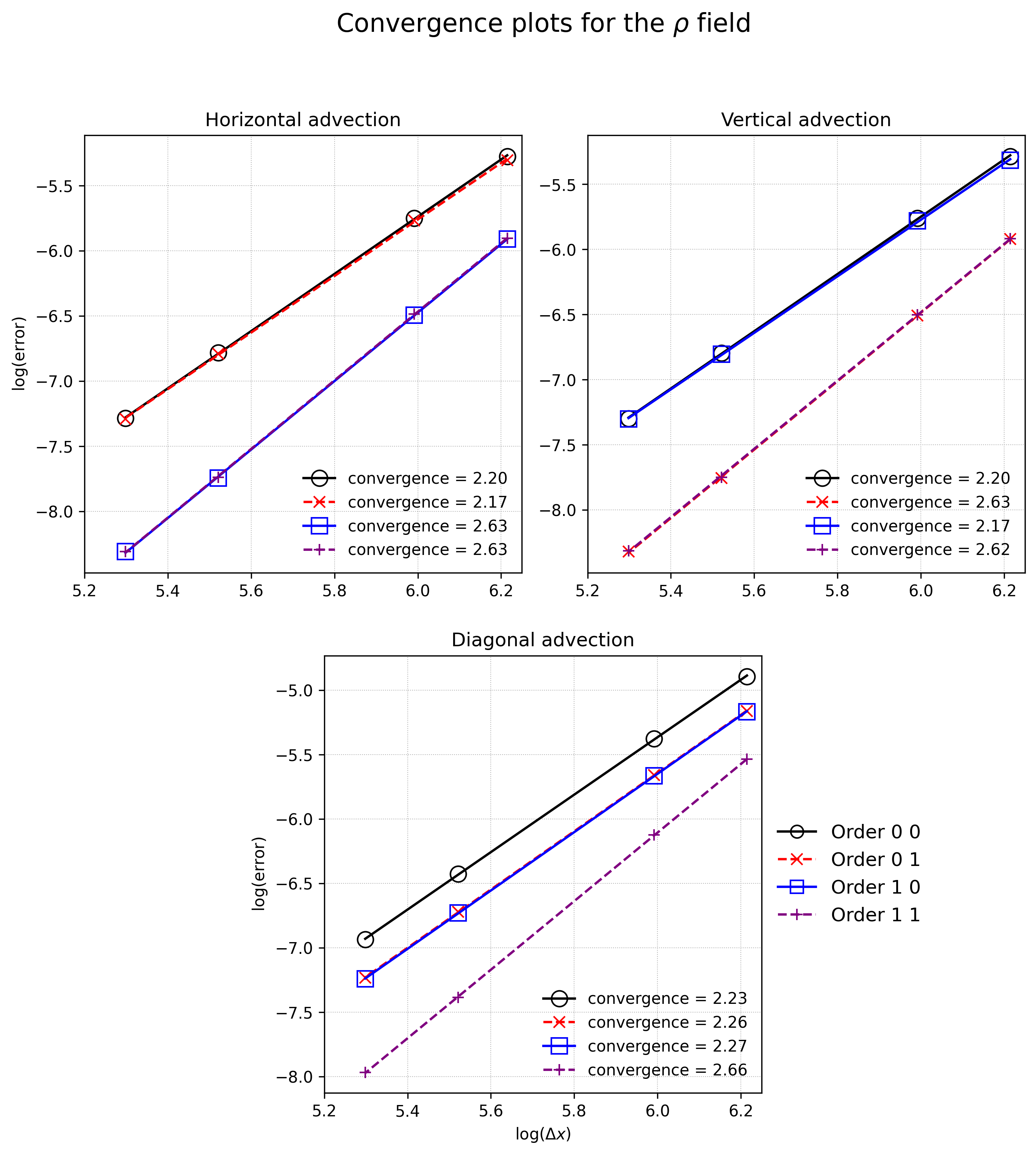}
    \caption{Convergence plots for each of the advection cases. Going left 
    to right the advecting velocities are, ($(100ms^{_1}, 0)$) for the horizontal,
    $(0, 100ms^{_1})$ for the vertical and, $(100ms^{_1}, 100ms^{_1})$ for the diagonal.}
    \label{fig:TV convergence}
\end{figure}
Figure~\ref{fig:TV convergence} demonstrates a clear relationship between convergence behaviour and the dominant direction of motion. For horizontal advection, the $(1,1)$ and $(1,0)$ cases converge more rapidly than the $(0,0)$ and $(0,1)$ cases. Conversely, under vertical advection, the vertical order governs convergence, with the corresponding pair performing better. When the vortex is advected diagonally across the domain, the $(0,0)$, $(1,0)$, and $(0,1)$ cases exhibit the same convergence rate but differ in their absolute error levels, while the $(1,1)$ case achieves both a lower error and a higher convergence rate. Increasing the order, therefore, systematically improves the convergence properties of the scheme for motion aligned with that direction.

\subsection{Sch\"ar Mountian wave} \label{sec: mountain wave}
So far, we have focused on the convergence order of the split-order scheme; however, accuracy considerations are equally important. As discussed in \cref{sec: dispersion}, increasing the polynomial order improves the representation of low wavenumbers, but introduces spectral gaps and degrades the representation of high wavenumbers. To assess these competing effects, we consider the test case of Sch"ar et al. \citep{scharNewTerrainFollowingVertical2002}, which excites a broad range of wavenumbers.
We initialise a balanced atmosphere over a mountain range expressed by,
\begin{equation}
    h(x) = h_0\exp\left(-\left(\frac{x}{a}\right)^2\right)\cos^2\left(\frac{\pi x}{\lambda}\right),
\end{equation}
where $h_0 = \qty{250}{\meter}$, $a = \qty{5e3}{\meter}$, and $\lambda = \qty{4e3}{\meter}$.
We calculate the simulation in a domain of dimensions $[-L/2 , L/2] \times [0, 30]\unit{\kilo\meter}$ where $L =\qty{100}{\kilo\meter}$.
The velocity is initially horizontal with $u=\qty{10}{\meter\per\second}$ and the potential temperature field has a profile of 
\begin{equation}
\theta_b = T_{surf}\exp\left(-\frac{N^2z}{g}\right), \quad T_{surf} = \qty{288}{\kelvin},
\end{equation}
where $N = \qty{0.01}{\per\second}$ and $g = \qty{9.810616}{\meter\per\square\second}$.
The atmosphere is initialised in hydrostatic balance following the procedure described in \citep{cotterCompatibleFiniteElement2022}.
In addition, we need to adjust the implicit-to-explicit ratio of our time-stepper by setting $\alpha=0.51$. This adjustment makes the time scheme slightly more implicit and dampens the most extreme grid-scale waves. This choice is standard for semi-implicit quasi-Newton time schemes implemented in full NWP models \cite{benacchioSemiimplicitSemiLagrangianModelling2016}.
\begin{figure}
    \centering
    \includegraphics[width=\linewidth]{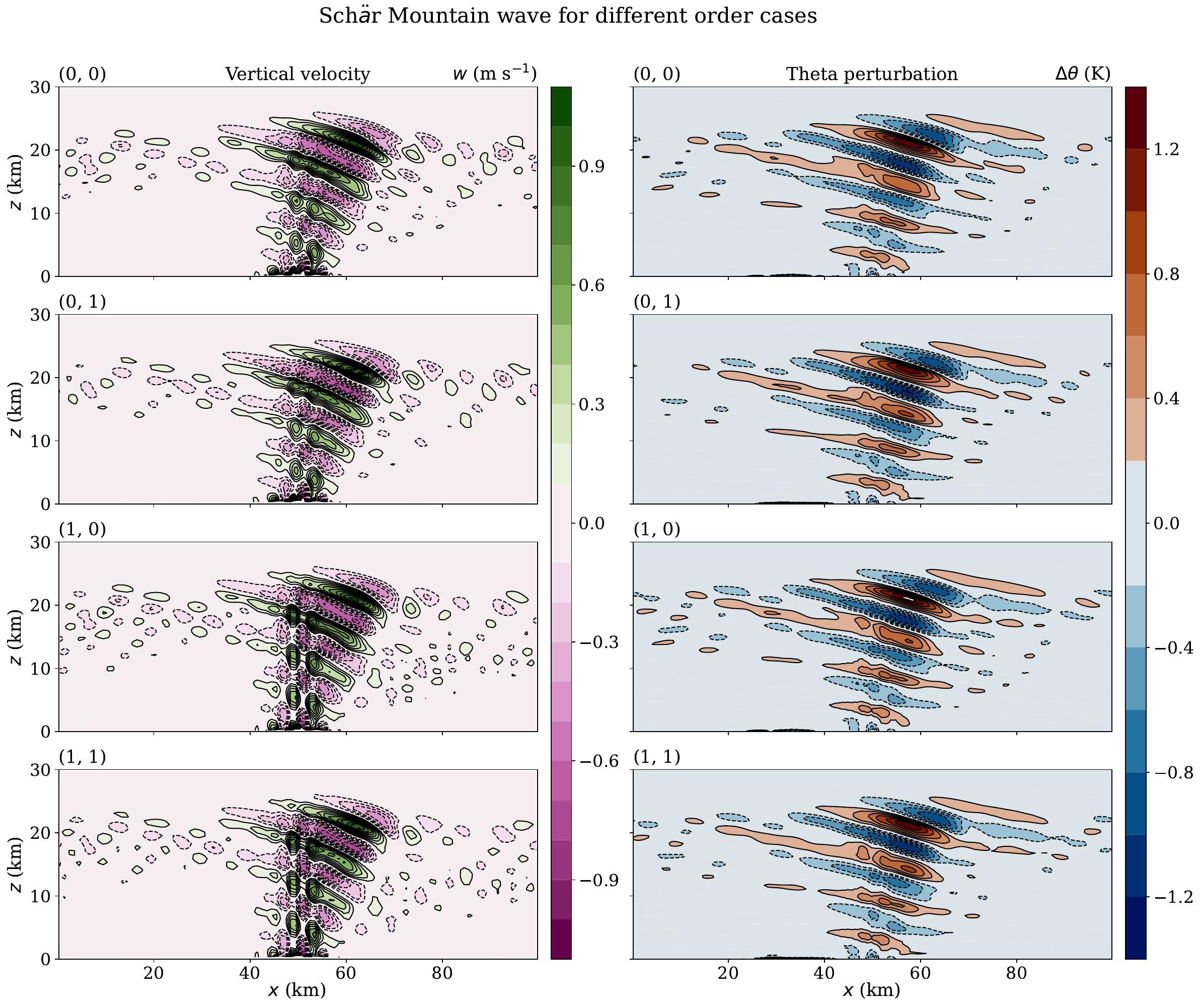}
    \caption{Plots of the vertical velocity (left) and theta perturbation
    (right) for the Sch\"ar mountain wave test case at $t=3600s$. Resolution for
    the four cases are, order $(1, 1)$: $(1000, 600)\unit{\meter}$, 
    order $(1, 0)$: $(1000, 300)\unit{\meter}$, order $(0, 1)$: $(500, 600)\unit{\meter}$
    and order $(0, 0)$: $(500, 300)\unit{\meter}$. 
    }
    \label{fig:schar mountain wave}
\end{figure}
We present results of the Sch\"ar mountain wave test in figure \ref{fig:schar mountain wave}.
Each case provides a solution that is consistent with those found in 
the literature
\citep{giraldoStudySpectralElement2008, cotterCompatibleFiniteElement2023, scharNewTerrainFollowingVertical2002}.
Whilst the theta perturbation is largely identical across the cases, for the 
vertical velocity, we can see that the higher horizontal cases $(1, 1)$ and 
$(1, 0)$ contain more pronounced vertical velocity
packets at $\qty{10}{\kilo\meter}$ and to the right of the mountain peak. 

\subsection{Baroclinic Wave} \label{sec: Baroclinic wave}
\begin{figure}
    \centering
    \includegraphics[width=\linewidth]{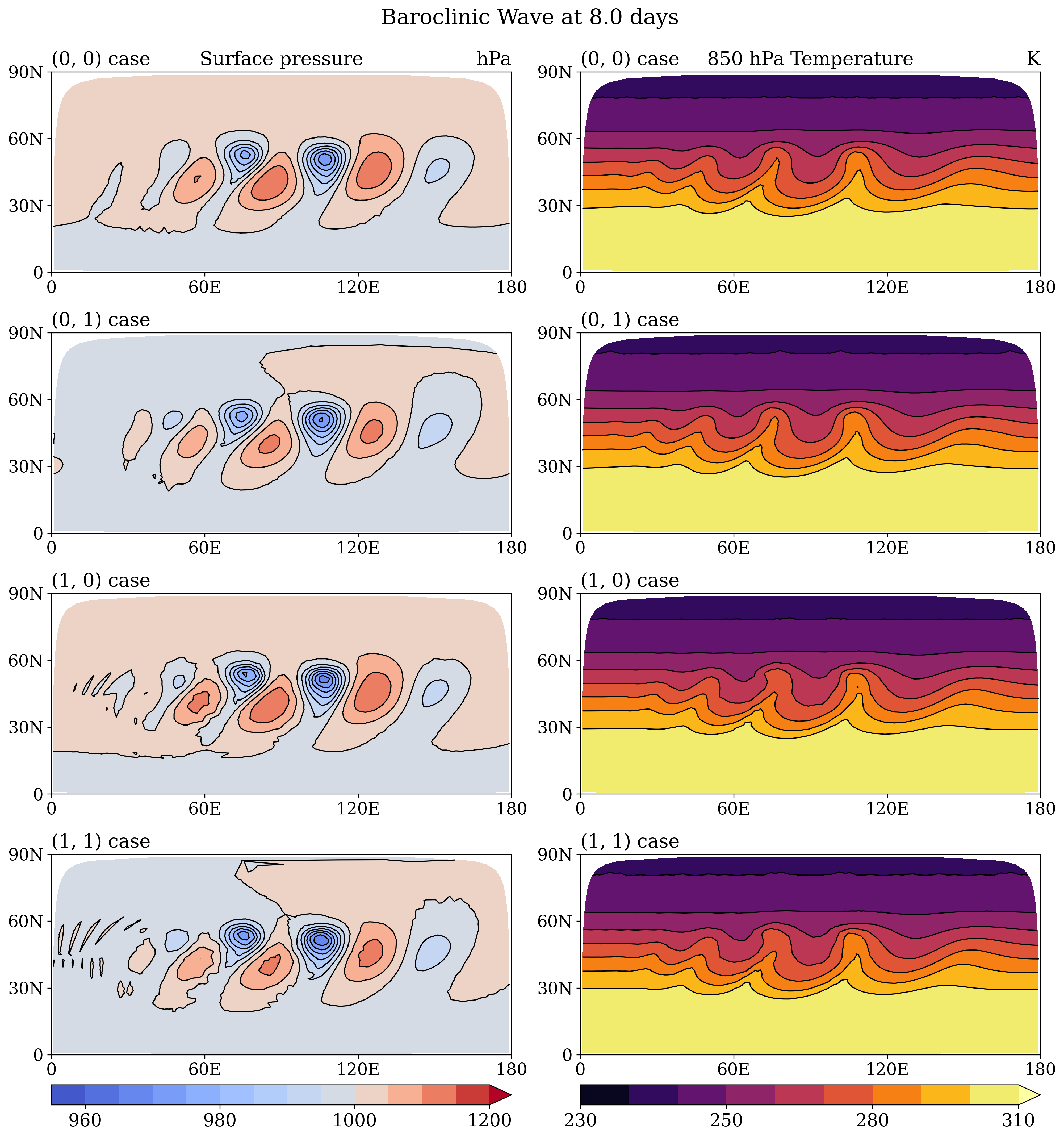}
    \caption{Plots of the surface pressure (left) and surface temperature (right) at 
    day 8 of integration. The rows correspond to the different finite element orders.}
    \label{fig:day 8 baroclinic}
\end{figure}

\begin{figure}
    \centering
    \includegraphics[width=\linewidth]{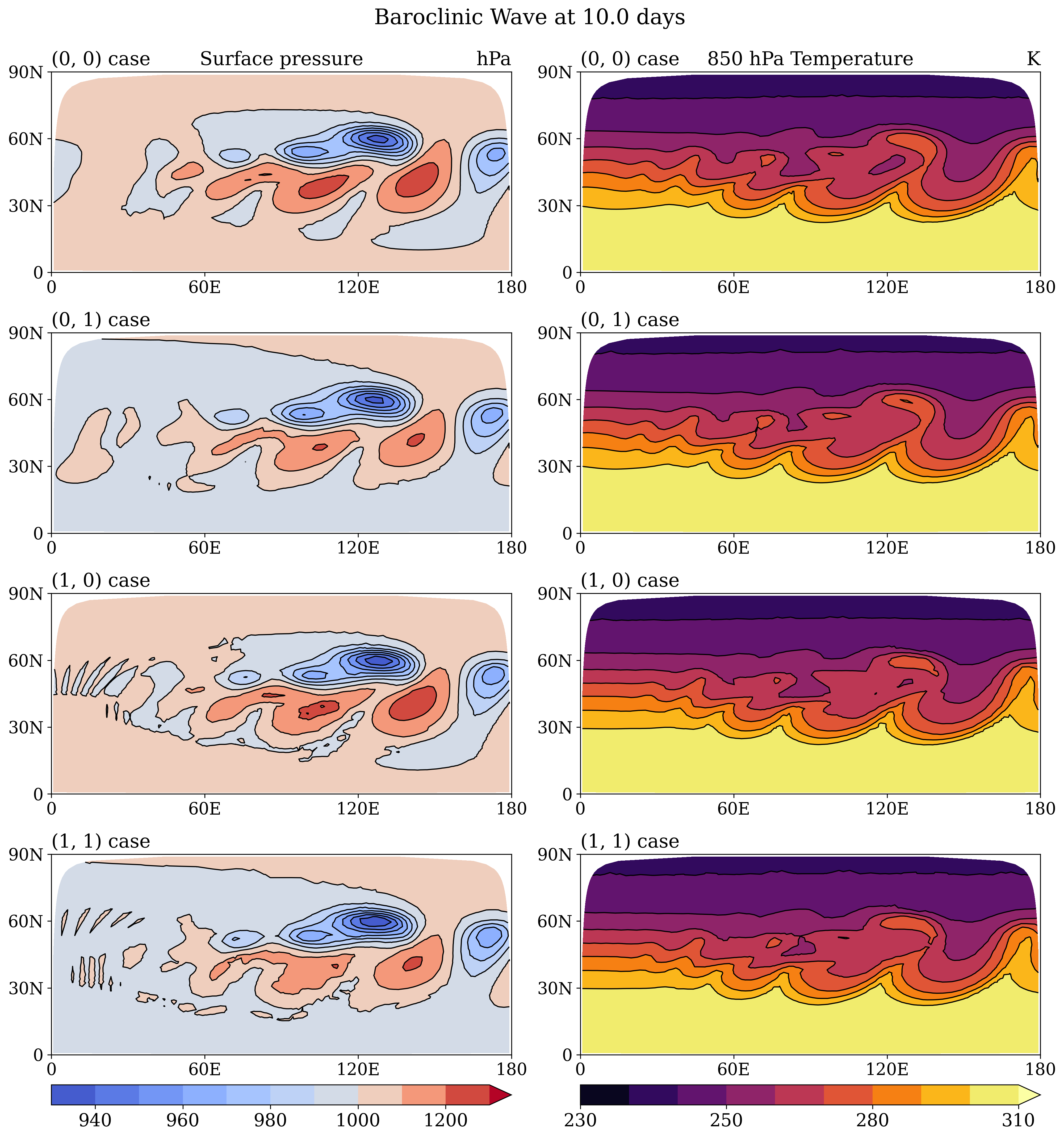}
    \caption{Plots of the surface pressure (left) and surface temperature (right) at 
    day 10 of integration. The rows correspond to the different finite element orders.}
    \label{fig:day 10 baroclinic}
\end{figure}

Finally, we demonstrate the successful implementation of the split-order schemes on an Earth-sized domain with a deep atmosphere using the dry baroclinic wave described in \cite{ullrichProposedBaroclinicWave2014}. We initialise the test on a cubed sphere grid with appropriate resolutions (see appendix \ref{A:resolutions}). Each simulation was run for 10 days, with a time step of $\Delta t = 450.0$s for all cases. Plots of the surface pressure and temperature at days 8 and ten 10 can be seen in figures 
\ref{fig:day 8 baroclinic} and \ref{fig:day 10 baroclinic}.

Each case produces a solution that agrees with the results present in the literature
\citep{ullrichProposedBaroclinicWave2014, woodInherentlyMassconservingSemiimplicit2014, melvinMixedFiniteelementFinitevolume2024, giraldoPerformanceStudyHorizontally2024, skamarockFullyCompressibleNonhydrostatic2021}
and the solutions exhibit pressure perturbations at the same locations.

\section{Conclusion}
The paper has focused on demonstrating the split-order scheme's accuracy and flexibility; however, the computational impacts of independently varying the horizontal and vertical polynomial orders remain to be fully characterised.
Through a discrete dispersion analysis of the compressible Boussinesq equations, we have examined how differing horizontal and vertical polynomial orders affect wave propagation. Higher-order discretisations reduce dispersion error at slow and intermediate wavenumbers but can introduce spectral gaps and increased error near the grid scale. The gravity wave mode was found to be particularly sensitive to changes in the vertical polynomial order, highlighting the importance of vertical discretisation choices in atmospheric models.
To further explore these findings, a suite of numerical test cases in section \ref{sec: test cases} demonstrated that increasing the polynomial order in a given direction improves convergence in that direction, with the horizontal order having a consistently larger impact under typical atmospheric conditions. Qualitative tests showed no evidence of spurious grid-scale wave activity associated with higher-order discretisations, despite concerns suggested by the dispersion analysis. Across all configurations, the split-order schemes reproduced results reported in the existing literature, supporting the use of polynomial order as a practical means of tuning the numerical behaviour of a dynamical core.
While the present study has focused primarily on accuracy considerations, it is important to note that computational efficiency and the interaction between split-order discretisations and physical parameterisations are also critical aspects of dynamical core design. Thus, the implications of split polynomial order for performance and physics–dynamics coupling remain open questions and represent important directions for future work.
\section*{Acknowledgements}
This research was funded by the Faculty of Environment Science and Economy at the University of Exeter. The authors would also like to thank Karina Kowalczyk for help with the initialisation of the travelling vortex test case.
For the purpose of open access, the authors have applied a ‘Creative Commons Attribution (CC BY) licence to any Author Accepted Manuscript version arising from this submission.

\appendix
\section{(0, 1) node location}\label{A:(0,1)locs}
The DoF locations for the finite elements with horizontal order 0 and vertical order 1 are shown in figure \ref{eqn:number_guide}.
\begin{figure}
    \centering
    \includegraphics[scale=1]{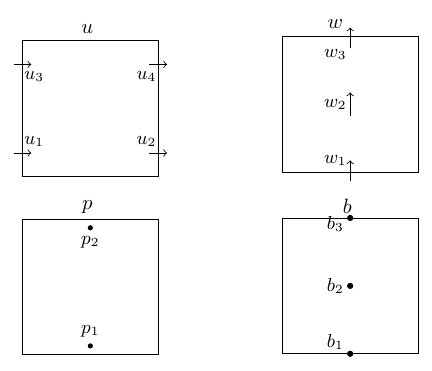}
    \caption{Locations of degrees of freedom for the (1,0) case for the 
             horizontal velocity field (u) shown in the top left followed by the 
             vertical velocity field (w) in the top right, pressure field (p) in 
             the bottom left and the buoyancy field (b) in the bottom right. The
             node numbering can be related to node location by using equation \eqref{eqn:number_guide}}
    \label{fig:DoF locations}
\end{figure}
\section{Dispersion matrices}\label{A:dispersion_matricies}
Details of the mass and discrete derivative matrices for the $(1,0)$ case can be found in section \ref{sec:1-0-matricies} and matrices for the $(0,1)$ case in section \ref{sec:0-1-matricies}. Details of the relevant matrices for the $(1,1)$ and $(0,0)$ cases can be found in \cite{melvinChoiceFunctionSpaces2018}.
\subsection{(1, 0) Case} \label{sec:1-0-matricies}
\begin{equation}\label{A:1-0-Mu-matrix}
    M_u =
    \frac{\Delta x \Delta z}{30}
    \left(
    \begin{array}{cc}
        \bigl(8 - 2\cos(k \Delta x)\bigr)\,\mathrm{e}^{\mathrm{i} l \Delta z / 2}
        &
        4\cos\!\left(\frac{k \Delta x}{2}\right)\,\mathrm{e}^{\mathrm{i} l \Delta z / 2}
        \\[0.5em]
        \bigl(2\,\mathrm{e}^{\mathrm{i} k \Delta x} + 2\bigr)\,\mathrm{e}^{\mathrm{i} l \Delta z / 2}
        &
        16\,\mathrm{e}^{\mathrm{i} k \Delta x / 2}\,\mathrm{e}^{\mathrm{i} l \Delta z / 2}
    \end{array}
    \right).
\end{equation}
\begin{equation}\label{A:1-0-MwMbQ-matrix-CP}
    M_w = M_b = Q =
    \frac{\Delta x \Delta z}{36}
    \left(
    \begin{array}{cc}
        4(\cos(l \Delta z) + 2)
        &
        2(\cos(l \Delta z) + 2)\,\mathrm{e}^{\mathrm{i} k \Delta x}
        \\[0.5em]
        2(\cos(l \Delta z) + 2)
        &
        4(\cos(l \Delta z) + 2)\,\mathrm{e}^{\mathrm{i} k \Delta x}
    \end{array}
    \right).
\end{equation}
\begin{equation}\label{A:1-0-Mp-matrix}
    M_p =
    \frac{\Delta x \Delta z}{6}
    \left(
    \begin{array}{cc}
        2\,\mathrm{e}^{\mathrm{i} l \Delta z / 2}
        &
        \mathrm{e}^{\mathrm{i} l \Delta z / 2}\,\mathrm{e}^{\mathrm{i} k \Delta x}
        \\[0.5em]
        \mathrm{e}^{\mathrm{i} l \Delta z / 2}
        &
        2\,\mathrm{e}^{\mathrm{i} l \Delta z / 2}\,\mathrm{e}^{\mathrm{i} k \Delta x}
    \end{array}
    \right).
\end{equation}
\begin{equation}\label{A:1-0-Dux-matrix-2}
    D^u_x =
    \frac{\Delta z}{6}
    \left(
    \begin{array}{cc}
        (5 - \mathrm{e}^{-\mathrm{i} k \Delta x})\,\mathrm{e}^{\mathrm{i} l \Delta z / 2}
        &
        (\mathrm{e}^{\mathrm{i} k \Delta x} - 5)\,\mathrm{e}^{\mathrm{i} l \Delta z / 2}
        \\[0.5em]
        -4\,\mathrm{e}^{\mathrm{i} l \Delta z / 2}
        &
        4\,\mathrm{e}^{\mathrm{i} l \Delta z / 2}\,\mathrm{e}^{\mathrm{i} k \Delta x}
    \end{array}
    \right).
\end{equation}

\begin{equation}\label{A:1-0-Dpx-matrix-2}
    D^p_x =
    \frac{\Delta z}{6}
    \left(
    \begin{array}{cc}
        (\mathrm{e}^{\mathrm{i} k \Delta x} - 5)\,\mathrm{e}^{\mathrm{i} l \Delta z / 2}
        &
        4\,\mathrm{e}^{\mathrm{i} l \Delta z / 2}\,\mathrm{e}^{\mathrm{i} k \Delta x / 2}
        \\[0.5em]
        (5\,\mathrm{e}^{\mathrm{i} k \Delta x} - 1)\,\mathrm{e}^{\mathrm{i} l \Delta z / 2}
        &
        -4\,\mathrm{e}^{\mathrm{i} l \Delta z / 2}\,\mathrm{e}^{\mathrm{i} k \Delta x / 2}
    \end{array}
    \right).
\end{equation}
\begin{equation}\label{A:1-0-Dwz-matrix-2}
    D^w_z =
    \frac{\Delta x}{6}
    \left(
    \begin{array}{cc}
        4\,\mathrm{i}\,\sin(l \Delta z / 2)
        &
        2\,\mathrm{i}\,\sin(l \Delta z / 2)\,\mathrm{e}^{\mathrm{i} k \Delta x}
        \\[0.5em]
        2\,\mathrm{i}\,\sin(l \Delta z / 2)
        &
        4\,\mathrm{i}\,\sin(l \Delta z / 2)\,\mathrm{e}^{\mathrm{i} k \Delta x}
    \end{array}
    \right).
\end{equation}
\begin{equation}\label{A:1-0-Dpz-matrix-2}
    D^p_z =
    \frac{\Delta x}{6}
    \left(
    \begin{array}{cc}
        2\,(\mathrm{e}^{\mathrm{i} l \Delta z} - 1)
        &
        (\mathrm{e}^{\mathrm{i} l \Delta z} - 1)\,\mathrm{e}^{\mathrm{i} k \Delta x}
        \\[0.5em]
        (\mathrm{e}^{\mathrm{i} l \Delta z} - 1)
        &
        2\,(\mathrm{e}^{\mathrm{i} l \Delta z} - 1)\,\mathrm{e}^{\mathrm{i} k \Delta x}
    \end{array}
    \right).
\end{equation}
\subsection{(0, 1) Case}\label{sec:0-1-matricies}
\begin{equation}\label{A:0-1-mass-u}
    M_u =
    \left(
    \begin{array}{cc}
        4\bigl(\cos(k\Delta x) + 2\bigr)
        &
        2\bigl(\cos(k\Delta x) + 2\bigr)\,\mathrm{e}^{\mathrm{i} l \Delta z}
        \\
        2\bigl(\cos(k\Delta x) + 2\bigr)
        &
        4\bigl(\cos(k\Delta x) + 2\bigr)\,\mathrm{e}^{\mathrm{i} l \Delta z}
    \end{array}
    \right).
\end{equation}
\begin{equation}
    M_w = M_b = Q =  
    \frac{\Delta x \Delta z}{30}
    \left(
    \begin{array}{cc}
        \bigl(8 - 2\cos(l\Delta z)\bigr)\,\mathrm{e}^{\mathrm{i} k \Delta x / 2}
        &
        4\cos\!\left(\tfrac{l\Delta z}{2}\right)\,\mathrm{e}^{\mathrm{i} k \Delta x / 2}
        \\[0.5em]
        \bigl(2\mathrm{e}^{\mathrm{i} l \Delta z} + 2\bigr)\,\mathrm{e}^{\mathrm{i} k \Delta x / 2}
        &
        16\,\mathrm{e}^{\mathrm{i} k \Delta x / 2}\,\mathrm{e}^{\mathrm{i} l \Delta z / 2}
    \end{array}
    \right).
\end{equation}
\begin{equation}\label{A:0-1-p-mass-matrix}
    M_p =
    \frac{\Delta x \Delta z}{6}
    \left(
    \begin{array}{cc}
        2\,\mathrm{e}^{\mathrm{i} k \Delta x / 2}
        &
        \mathrm{e}^{\mathrm{i} k \Delta x / 2}\,\mathrm{e}^{\mathrm{i} l \Delta z}
        \\[0.5em]
        \mathrm{e}^{\mathrm{i} k \Delta x / 2}
        &
        2\,\mathrm{e}^{\mathrm{i} k \Delta x / 2}\,\mathrm{e}^{\mathrm{i} l \Delta z}
    \end{array}
    \right).
\end{equation}

\begin{equation}\label{A:0-1-Dux-matrix}
    D^u_x =
    \frac{\Delta z}{6}
    \left(
    \begin{array}{cc}
        4\,\mathrm{i}\,\sin\!\left(\frac{k \Delta x}{2}\right)
        &
        2\,\mathrm{i}\,\sin\!\left(\frac{k \Delta x}{2}\right)\,\mathrm{e}^{\mathrm{i} l \Delta z}
        \\[0.5em]
        2\,\mathrm{i}\,\sin\!\left(\frac{k \Delta x}{2}\right)
        &
        4\,\mathrm{i}\,\sin\!\left(\frac{k \Delta x}{2}\right)\,\mathrm{e}^{\mathrm{i} l \Delta z}
    \end{array}
    \right).
\end{equation}

\begin{equation}\label{A:Dpx-matrix}
    D^p_x =
    \frac{\Delta z}{6}
    \left(
    \begin{array}{cc}
        2\left(\mathrm{e}^{\mathrm{i} k \Delta x} - 1\right)
        &
        \left(\mathrm{e}^{\mathrm{i} k \Delta x} - 1\right)\,\mathrm{e}^{\mathrm{i} l \Delta z}
        \\[0.5em]
        \left(\mathrm{e}^{\mathrm{i} k \Delta x} - 1\right)
        &
        2\left(\mathrm{e}^{\mathrm{i} k \Delta x} - 1\right)\,\mathrm{e}^{\mathrm{i} l \Delta z}
    \end{array}
    \right).
\end{equation}
\begin{equation}\label{A:0-1-Dwz-matrix}
    D^w_z =
    \frac{\Delta x}{6}
    \left(
    \begin{array}{cc}
        \bigl(5 - \mathrm{e}^{-\mathrm{i} l \Delta z}\bigr)\,\mathrm{e}^{\mathrm{i} k \Delta x / 2}
        &
        \bigl(\mathrm{e}^{\mathrm{i} l \Delta z} - 5\bigr)\,\mathrm{e}^{\mathrm{i} k \Delta x / 2}
        \\[0.5em]
        -4\,\mathrm{e}^{\mathrm{i} k \Delta x / 2}
        &
        4\,\mathrm{e}^{\mathrm{i} k \Delta x / 2}\,\mathrm{e}^{\mathrm{i} l \Delta z}
    \end{array}
    \right).
\end{equation}
\begin{equation}\label{A:Dpz-matrix}
    D^p_z =
    \frac{\Delta x}{6}
    \left(
    \begin{array}{cc}
        \bigl(\mathrm{e}^{\mathrm{i} l \Delta z} - 5\bigr)\,\mathrm{e}^{\mathrm{i} k \Delta x / 2}
        &
        4\,\mathrm{e}^{\mathrm{i} k \Delta x / 2}\,\mathrm{e}^{\mathrm{i} l \Delta z / 2}
        \\[0.5em]
        \bigl(5\,\mathrm{e}^{\mathrm{i} l \Delta z} - 1\bigr)\,\mathrm{e}^{\mathrm{i} k \Delta x / 2}
        &
        -4\,\mathrm{e}^{\mathrm{i} k \Delta x / 2}\,\mathrm{e}^{\mathrm{i} l \Delta z / 2}
    \end{array}
    \right).
\end{equation}
\section{Acoustic mode statistics} \label{A:accoustic_stats}
Details of the error statistics for the acoustic mode can be found in table \eqref{tab: acoustic error statistics}.
\begin{table}
\centering
\caption{Global error measured for the acoustic mode, each case has the normalised 
$L_2$ error is calculated, as well as the maximum and minimum error. The higher-order cases 
perform better overall, however, have larger maximum errors.}
{\renewcommand{\arraystretch}{2}%
\begin{tabular}{lrrr}
    \hline
     case   &   Normalised $L_2$ error norm &   Max error &   Min Error \\
    \hline
     (0, 0) &                0.000482784 &   0.0404953 & 2.08984e-09 \\
     (0, 1) &                0.000465869 &   0.0824462 & 3.2898e-10  \\
     (1, 0) &                0.000465834 &   0.0824227 & 9.60074e-10 \\
     (1, 1) &                0.000400848 &   0.0826268 & 8.27239e-11 \\
    \hline
\end{tabular}}
\label{tab: acoustic error statistics}
\end{table}
\section{Timestep pseudocode} \label{A:pseudo-code}
Pseudocode for the timestep procedure is found in \cref{alg:timestep}
\begin{algorithm} \label{alg:timestep}
    \caption{Pseudocode of a single Semi-Implicit Quasi-Newton time step, the initial state is given as 
    $X^0$ with $k_{\max}$ and $I_{\max}$ being the number of outer and inner loops, respectively.} 
    \label{psudo:timestep1}
    \begin{algorithmic}[1]
    \State Set: $X^{n+1} \gets X^n$
    \State Explicit forcing: $X_f \gets X^n + (1-\alpha)\Delta tF(X^n)$ 
    \State $k \gets 0$
    \State \textbf{Outer loop}
    \While{$k < k_{\max}$}:    
        \State Update advecting velocity: $\bar{\mathbf{u}} \gets \frac{1}{2}(\mathbf{u}^{n+1} + \mathbf{u}^n)$
        \State Advect state: $X_a \gets A_{\bar{\mathbf{u}}}(X_f)$
        \State $i \gets 0$
        \State \textbf{Inner loop}
        \While{$i < I_{\max}$} 
            \State Implicit forcing: $X_f \gets X_a + \alpha \Delta tF(X^n)$
            \State Calculate residual: $X_r \gets X_f - X^{n+1}$
            \State Solve: $S(\Delta X) = X_r$ for $\Delta X$
            \State Increment $X^{n+1} \gets X^{n+1} + \Delta X$
            \State $i \gets i+1$
        \EndWhile
        \State $k \gets k+1$
    \EndWhile
    \State Advance time step: $X^n \gets X^{n+1}$
    \end{algorithmic}
\end{algorithm}
\section{Test case resolutions} \label{A:resolutions}
As mentioned in section \ref{sec: test cases}, we ensure that in each test case the average distance between DoFs is roughly kept the same. Table \ref{tab:res-scaling} details the horizontal and vertical resolutions for each test case and each order.
\begin{table}[]
\centering
\begin{tabular}{@{}lllll@{}}
\toprule
Test Case         & (0, 0)         & (0, 1)         & (1, 0)         & (1, 1)         \\ \midrule
Gravity wave      & (1000m, 1000m) & (1000m, 2000m) & (2000m, 1000m) & (2000m, 2000m) \\
Travelling vortex & (100m, 100m)   & (100m, 200m)   & (200m, 100m)   & (200m, 200m)   \\
Sch\"ar Wave        & (500m, 300m)   & (500m, 600m)   & (1000,m 300m)  & (1000m, 600m)  \\ \bottomrule
\end{tabular}
\caption{Reference table for resolutions of idealised test cases for the different order cases}
\label{tab:res-scaling}
\end{table}
\bibliographystyle{wileyqj}
\bibliography{SplitOrderPaper-bibliography}

@article{adamsLFRicMeetingChallenges2019,
  title = {{{LFRic}}: {{Meeting}} the Challenges of Scalability and Performance Portability in {{Weather}} and {{Climate}} Models},
  author = {Adams, S. V. and Ford, R. W. and Hambley, M. and Hobson, J. M. and Kav{\v c}i{\v c}, I. and Maynard, C. M. and Melvin, T. and M{\"u}ller, E. H. and Mullerworth, S. and Porter, A. R. and Rezny, M. and Shipway, B. J. and Wong, R.},
  year = 2019,
  journal = {Journal of Parallel and Distributed Computing},
  volume = {132},
  pages = {383--396},
  issn = {0743-7315},
  doi = {10.1016/j.jpdc.2019.02.007}
}

@incollection{arakawaComputationalDesignBasic1977,
  title = {Computational {{Design}} of the {{Basic Dynamical Processes}} of the {{UCLA General Circulation Model}}},
  booktitle = {Methods in {{Computational Physics}}: {{Advances}} in {{Research}} and {{Applications}}},
  author = {Arakawa, {\relax AKIO} and Lamb, VIVIAN R.},
  editor = {Chang, {\relax JULIUS}},
  year = 1977,
  month = jan,
  series = {General {{Circulation Models}} of the {{Atmosphere}}},
  volume = {17},
  pages = {173--265},
  publisher = {Elsevier},
  doi = {10.1016/B978-0-12-460817-7.50009-4},
  langid = {english}
}

@incollection{arnoldDifferentialComplexesStability2006,
  title = {Differential {{Complexes}} and {{Stability}} of {{Finite Element Methods I}}. {{The}} de {{Rham Complex}}},
  booktitle = {Compatible {{Spatial Discretizations}}},
  author = {Arnold, Douglas N. and Falk, Richard S. and Winther, Ragnar},
  editor = {Arnold, Douglas N. and Bochev, Pavel B. and Lehoucq, Richard B. and Nicolaides, Roy A. and Shashkov, Mikhail},
  year = 2006,
  volume = {142},
  pages = {23--46},
  publisher = {Springer New York},
  address = {New York, NY},
  doi = {10.1007/0-387-38034-5_2},
  isbn = {978-0-387-30916-3},
  langid = {english}
}

@article{bendallImprovingAccuracyDiscretisations2023,
  title = {Improving the Accuracy of Discretisations of the Vector Transport Equation on the Lowest-Order Quadrilateral {{Raviart-Thomas}} Finite Elements},
  author = {Bendall, T. M. and Wimmer, G. A.},
  year = 2023,
  month = feb,
  journal = {Journal of Computational Physics},
  volume = {474},
  pages = {111834},
  issn = {0021-9991},
  doi = {10.1016/j.jcp.2022.111834}
}

@article{bendallRecoveredSpaceAdvection2019,
  title = {The {{Recovered Space Advection Scheme}} for {{Lowest-Order Compatible Finite Element Methods}}},
  author = {Bendall, Thomas M. and Cotter, Colin J. and Shipton, Jemma},
  year = 2019,
  month = aug,
  journal = {Journal of Computational Physics},
  volume = {390},
  eprint = {1811.06956},
  primaryclass = {math},
  pages = {342--358},
  issn = {00219991},
  doi = {10.1016/j.jcp.2019.04.013},
  archiveprefix = {arXiv}
}

@article{bendallSolutionTrilemmaMoist2023,
  title = {A Solution to the Trilemma of the Moist {{Charney}}--{{Phillips}} Staggering},
  author = {Bendall, Thomas M. and Wood, Nigel and Thuburn, John and Cotter, Colin J.},
  year = 2023,
  journal = {Quarterly Journal of the Royal Meteorological Society},
  volume = {149},
  number = {750},
  pages = {262--276},
  issn = {1477-870X},
  doi = {10.1002/qj.4406},
  copyright = {\copyright{} 2022 Crown copyright. Quarterly Journal of the Royal Meteorological Society \copyright{} 2022 Royal Meteorological Society. This article is published with the permission of the Controller of HMSO and the King's Printer for Scotland.},
  langid = {english}
}

@article{charneyNUMERICALINTEGRATIONQUASIGEOSTROPHIC1953,
  title = {{{NUMERICAL INTEGRATION OF THE QUASI-GEOSTROPHIC EQUATIONS FOR BAROTROPIC AND SIMPLE BAROCLINIC FLOWS}}},
  author = {Charney, J. G. and Phillips, N. A.},
  year = 1953,
  month = apr,
  journal = {Journal of the Atmospheric Sciences},
  volume = {10},
  number = {2},
  pages = {71--99},
  publisher = {American Meteorological Society},
  issn = {1520-0469},
  doi = {10.1175/1520-0469(1953)010<0071:NIOTQG>2.0.CO;2},
  chapter = {Journal of the Atmospheric Sciences},
  langid = {english}
}

@misc{cotterCompatibleFiniteElement2022,
  title = {A Compatible Finite Element Discretisation for the Nonhydrostatic Vertical Slice Equations},
  author = {Cotter, C. J. and Shipton, J.},
  year = 2022,
  month = oct,
  eprint = {2210.07861},
  primaryclass = {cs, math},
  doi = {10.48550/arXiv.2210.07861},
  archiveprefix = {arXiv},
  howpublished = {http://arxiv.org/abs/2210.07861},
  pubstate = {prepublished}
}

@misc{cotterCompatibleFiniteElement2023,
  title = {Compatible Finite Element Methods for Geophysical Fluid Dynamics},
  author = {Cotter, Colin J.},
  year = 2023,
  month = mar,
  number = {arXiv:2302.13337},
  eprint = {2302.13337},
  primaryclass = {physics},
  publisher = {arXiv},
  doi = {10.48550/arXiv.2302.13337},
  archiveprefix = {arXiv}
}

@article{cotterFiniteElementExterior2014,
  title = {A Finite Element Exterior Calculus Framework for the Rotating Shallow-Water Equations},
  author = {Cotter, C. J. and Thuburn, J.},
  year = 2014,
  month = jan,
  journal = {1526},
  publisher = {ACADEMIC PRESS INC ELSEVIER SCIENCE},
  issn = {0021-9991},
  doi = {10.1016/j.jcp.2013.10.008},
  langid = {english}
}

@article{cotterMixedFiniteElements2012,
  title = {Mixed Finite Elements for Numerical Weather Prediction},
  author = {Cotter, C. J. and Shipton, J.},
  year = 2012,
  month = aug,
  journal = {Journal of Computational Physics},
  volume = {231},
  number = {21},
  pages = {7076--7091},
  issn = {0021-9991},
  doi = {10.1016/j.jcp.2012.05.020},
  langid = {english}
}

@book{durranNumericalMethodsFluid2010,
  title = {Numerical {{Methods}} for {{Fluid Dynamics}}},
  author = {Durran, Dale R.},
  year = 2010,
  series = {Texts in {{Applied Mathematics}}},
  volume = {32},
  publisher = {Springer},
  address = {New York, NY},
  doi = {10.1007/978-1-4419-6412-0},
  isbn = {978-1-4419-6411-3 978-1-4419-6412-0}
}

@article{eldredDispersionAnalysisCompatible2018,
  title = {Dispersion Analysis of Compatible {{Galerkin}} Schemes for the {{1D}} Shallow Water Model},
  author = {Eldred, Christopher and Le Roux, Daniel Y.},
  year = 2018,
  month = oct,
  journal = {Journal of Computational Physics},
  volume = {371},
  pages = {779--800},
  issn = {0021-9991},
  doi = {10.1016/j.jcp.2018.06.007}
}

@article{giraldoPerformanceStudyHorizontally2024,
  title = {A Performance Study of Horizontally Explicit Vertically Implicit ({{HEVI}}) Time-Integrators for Non-Hydrostatic Atmospheric Models},
  author = {Giraldo, Francis X. and {de Bragan{\c c}a Alves}, Felipe A. V. and Kelly, James F. and Kang, Soonpil and Reinecke, P. Alex},
  year = 2024,
  month = oct,
  journal = {Journal of Computational Physics},
  volume = {515},
  pages = {113275},
  issn = {0021-9991},
  doi = {10.1016/j.jcp.2024.113275}
}

@article{giraldoStudySpectralElement2008,
  title = {A Study of Spectral Element and Discontinuous {{Galerkin}} Methods for the {{Navier}}--{{Stokes}} Equations in Nonhydrostatic Mesoscale Atmospheric Modeling: {{Equation}} Sets and Test Cases},
  shorttitle = {A Study of Spectral Element and Discontinuous {{Galerkin}} Methods for the {{Navier}}--{{Stokes}} Equations in Nonhydrostatic Mesoscale Atmospheric Modeling},
  author = {Giraldo, F. X. and Restelli, M.},
  year = 2008,
  month = apr,
  journal = {Journal of Computational Physics},
  volume = {227},
  number = {8},
  pages = {3849--3877},
  issn = {0021-9991},
  doi = {10.1016/j.jcp.2007.12.009}
}

@article{guerraHighorderStaggeredFiniteelement2016,
  title = {A High-Order Staggered Finite-Element Vertical Discretization for Non-Hydrostatic Atmospheric Models},
  author = {Guerra, Jorge E. and Ullrich, Paul A.},
  year = 2016,
  month = jun,
  journal = {Geoscientific Model Development},
  volume = {9},
  number = {5},
  pages = {2007--2029},
  publisher = {Copernicus GmbH},
  issn = {1991-959X},
  doi = {10.5194/gmd-9-2007-2016},
  langid = {english}
}

@article{hannahSeparatingPhysicsDynamics2021,
  title = {Separating {{Physics}} and {{Dynamics Grids}} for {{Improved Computational Efficiency}} in {{Spectral Element Earth System Models}}},
  author = {Hannah, Walter M. and Bradley, Andrew M. and Guba, Oksana and Tang, Qi and Golaz, Jean-Christophe and Wolfe, Jon},
  year = 2021,
  journal = {Journal of Advances in Modeling Earth Systems},
  volume = {13},
  number = {7},
  pages = {e2020MS002419},
  issn = {1942-2466},
  doi = {10.1029/2020MS002419},
  copyright = {\copyright{} 2021. The Authors. Journal of Advances in Modeling Earth Systems published by Wiley Periodicals LLC on behalf of American Geophysical Union.},
  langid = {english}
}

@misc{hartneyExploringFormsMoist2025,
  title = {Exploring Forms of the Moist Shallow Water Equations Using a New Compatible Finite Element Discretisation},
  author = {Hartney, Nell and Bendall, Thomas M. and Shipton, Jemma},
  year = 2025,
  month = may,
  eprint = {2409.07182},
  primaryclass = {math},
  doi = {10.48550/arXiv.2409.07182},
  archiveprefix = {arXiv},
  howpublished = {http://arxiv.org/abs/2409.07182},
  pubstate = {prepublished}
}

@article{herringtonPhysicsDynamicsCoupling2019,
  title = {Physics--{{Dynamics Coupling}} with {{Element-Based High-Order Galerkin Methods}}: {{Quasi-Equal-Area Physics Grid}}},
  shorttitle = {Physics--{{Dynamics Coupling}} with {{Element-Based High-Order Galerkin Methods}}},
  author = {Herrington, Adam R. and Lauritzen, Peter H. and Taylor, Mark A. and Goldhaber, Steve and Eaton, Brian E. and Bacmeister, Julio T. and Reed, Kevin A. and Ullrich, Paul A.},
  journal = "Monthly Weather Review",
  year = 2019,
  month = jan,
  doi = {10.1175/MWR-D-18-0136.1},
  chapter = {Monthly Weather Review},
  langid = {english}
}

@article{kadiogluFourthorderAuxiliaryVariable2008,
  title = {A Fourth-Order Auxiliary Variable Projection Method for Zero-{{Mach}} Number Gas Dynamics},
  author = {Kadioglu, Samet Y. and Klein, Rupert and Minion, Michael L.},
  year = 2008,
  month = jan,
  journal = {Journal of Computational Physics},
  volume = {227},
  pages = {2012--2043},
  issn = {0021-9991},
  doi = {10.1016/j.jcp.2007.10.008}
}

@article{lerouxAnalysisNumericallyInduced2007,
  title = {Analysis of {{Numerically Induced Oscillations}} in {{2D Finite}}-{{Element Shallow}}-{{Water Models Part I}}: {{Inertia}}-{{Gravity Waves}}},
  shorttitle = {Analysis of {{Numerically Induced Oscillations}} in {{2D Finite}}-{{Element Shallow}}-{{Water Models Part I}}},
  author = {Le Roux, Daniel Y. and Rostand, Virgile and Pouliot, Benoit},
  year = 2007,
  month = jan,
  journal = {SIAM Journal on Scientific Computing},
  volume = {29},
  number = {1},
  pages = {331--360},
  publisher = {{Society for Industrial and Applied Mathematics}},
  issn = {1064-8275},
  doi = {10.1137/060650106}
}

@article{melvinChoiceFunctionSpaces2018,
  title = {Choice of Function Spaces for Thermodynamic Variables in Mixed Finite Element Methods},
  author = {Melvin, Thomas and Benacchio, Tommaso and Thuburn, John and Cotter, Colin},
  year = 2018,
  month = aug,
  journal = {Quarterly Journal of the Royal Meteorological Society},
  volume = {144},
  pages = {900--916},
  doi = {10.1002/qj.3268}
}

@article{melvinInherentlyMassconservingIterative2010,
  title = {An Inherently Mass-Conserving Iterative Semi-Implicit Semi-{{Lagrangian}} Discretization of the Non-Hydrostatic Vertical-Slice Equations},
  author = {Melvin, Thomas and Dubal, Mark and Wood, Nigel and Staniforth, Andrew and Zerroukat, Mohamed},
  year = 2010,
  journal = {Quarterly Journal of the Royal Meteorological Society},
  volume = {136},
  number = {648},
  pages = {799--814},
  issn = {1477-870X},
  doi = {10.1002/qj.603},
  langid = {english}
}

@article{melvinMixedFiniteelementFinitevolume2019,
  title = {A Mixed Finite-Element, Finite-Volume, Semi-Implicit Discretisation for Atmospheric Dynamics: {{Cartesian}} Geometry},
  shorttitle = {A Mixed Finite-Element, Finite-Volume, Semi-Implicit Discretisation for Atmospheric Dynamics},
  author = {Melvin, Thomas and Benacchio, Tommaso and Shipway, Ben and Wood, Nigel and Thuburn, John and Cotter, Colin},
  year = 2019,
  month = oct,
  journal = {Quarterly Journal of the Royal Meteorological Society},
  volume = {145},
  pages = {2835--2853},
  doi = {10.1002/qj.3501}
}

@article{melvinMixedFiniteelementFinitevolume2024,
  title = {A Mixed Finite-Element, Finite-Volume, Semi-Implicit Discretisation for Atmospheric Dynamics: {{Spherical}} Geometry},
  shorttitle = {A Mixed Finite-Element, Finite-Volume, Semi-Implicit Discretisation for Atmospheric Dynamics},
  author = {Melvin, Thomas and Shipway, Ben and Wood, Nigel and Benacchio, Tommaso and Bendall, Thomas and Boutle, Ian and Brown, Alex and Johnson, Christine and Kent, James and Pring, Stephen and Smith, Chris and Zerroukat, Mohamed and Cotter, Colin and Thuburn, John},
  year = 2024,
  journal = {Quarterly Journal of the Royal Meteorological Society},
  volume = {150},
  number = {764},
  pages = {4252--4269},
  issn = {1477-870X},
  doi = {10.1002/qj.4814},
  langid = {english}
}

@article{melvinTwodimensionalMixedFiniteelement2014,
  title = {A Two-Dimensional Mixed Finite-Element Pair on Rectangles},
  author = {Melvin, Thomas and Staniforth, Andrew and Cotter, Colin},
  year = 2014,
  journal = {Quarterly Journal of the Royal Meteorological Society},
  volume = {140},
  number = {680},
  pages = {930--942},
  issn = {1477-870X},
  doi = {10.1002/qj.2189},
  langid = {english}
}

@article{nairComputationalAspectsScalable2009,
  title = {Computational Aspects of a Scalable High-Order Discontinuous {{Galerkin}} Atmospheric Dynamical Core},
  author = {Nair, R. D. and Choi, H. -W. and Tufo, H. M.},
  year = 2009,
  month = feb,
  journal = {Computers \& Fluids},
  volume = {38},
  number = {2},
  pages = {309--319},
  issn = {0045-7930},
  doi = {10.1016/j.compfluid.2008.04.006}
}

@article{nataleCompatibleFiniteElement2016,
  title = {Compatible Finite Element Spaces for Geophysical Fluid Dynamics},
  author = {Natale, Andrea and Shipton, Jemma and Cotter, Colin J},
  year = 2016,
  month = jan,
  journal = {Dynamics and Statistics of the Climate System},
  volume = {1},
  number = {1},
  pages = {dzw005},
  issn = {2059-6987},
  doi = {10.1093/climsys/dzw005}
}

@article{scharNewTerrainFollowingVertical2002,
  title = {A {{New Terrain-Following Vertical Coordinate Formulation}} for {{Atmospheric Prediction Models}}},
  author = {Sch{\"a}r, Christoph and Leuenberger, Daniel and Fuhrer, Oliver and L{\"u}thi, Daniel and Girard, Claude},
  year = 2002,
  month = oct,
  journal = {American Meteorological Society},
  volume = {130},
  number = {10},
  pages = {2459--2480},
  issn = {1520-0493},
  chapter = {Monthly Weather Review},
  langid = {english}
}

@article{schmidtEffectsVerticalGrid2024,
  title = {Effects of Vertical Grid Spacing on the Climate Simulated in the {{ICON-Sapphire}} Global Storm-Resolving Model},
  author = {Schmidt, Hauke and Rast, Sebastian and Bao, Jiawei and Cassim, Amrit and Fang, Shih-Wei and {Jimenez-de la Cuesta}, Diego and Keil, Paul and Kluft, Lukas and Kroll, Clarissa and Lang, Theresa and Niemeier, Ulrike and Schneidereit, Andrea and Williams, Andrew I. L. and Stevens, Bjorn},
  year = 2024,
  month = feb,
  journal = {Geoscientific Model Development},
  volume = {17},
  number = {4},
  pages = {1563--1584},
  publisher = {Copernicus GmbH},
  issn = {1991-959X},
  doi = {10.5194/gmd-17-1563-2024},
  langid = {english}
}

@article{skamarockEfficiencyAccuracyKlempWilhelmson1994,
  title = {Efficiency and {{Accuracy}} of the {{Klemp-Wilhelmson Time-Splitting Technique}}},
  author = {Skamarock, William C. and Klemp, Joseph B.},
  year = 1994,
  month = nov,
  journal = {Monthly Weather Review},
  volume = {122},
  number = {11},
  pages = {2623--2630},
  publisher = {American Meteorological Society},
  issn = {1520-0493, 0027-0644},
  doi = {10.1175/1520-0493(1994)122<2623:EAAOTK>2.0.CO;2},
  chapter = {Monthly Weather Review},
  langid = {english}
}

@article{staniforthHorizontalGridsGlobal2012,
  title = {Horizontal Grids for Global Weather and Climate Prediction Models: A Review: {{Horizontal Grids}} for {{Global Models}}: A {{Review}}},
  shorttitle = {Horizontal Grids for Global Weather and Climate Prediction Models},
  author = {Staniforth, Andrew and Thuburn, John},
  year = 2012,
  month = jan,
  journal = {Quarterly Journal of the Royal Meteorological Society},
  volume = {138},
  number = {662},
  pages = {1--26},
  issn = {00359009},
  doi = {10.1002/qj.958},
  langid = {english}
}

@article{ullrichProposedBaroclinicWave2014,
  title = {A Proposed Baroclinic Wave Test Case for Deep- and Shallow-Atmosphere Dynamical Cores},
  author = {Ullrich, Paul A. and Melvin, Thomas and Jablonowski, Christiane and Staniforth, Andrew},
  year = 2014,
  journal = {Quarterly Journal of the Royal Meteorological Society},
  volume = {140},
  number = {682},
  pages = {1590--1602},
  issn = {1477-870X},
  doi = {10.1002/qj.2241},
  copyright = {\copyright{} 2013 Royal Meteorological Society and Crown Copyright, the Met Office. Quarterly Journal of the Royal Meteorological Society \copyright{} 2013 Royal Meteorological Society},
  langid = {english}
}

@article{woodInherentlyMassconservingSemiimplicit2014,
  title = {An Inherently Mass-Conserving Semi-Implicit Semi-{{Lagrangian}} Discretization of the Deep-Atmosphere Global Non-Hydrostatic Equations},
  author = {Wood, Nigel and Staniforth, Andrew and White, Andy and Allen, Thomas and Diamantakis, Michail and Gross, Markus and Melvin, Thomas and Smith, Chris and Vosper, Simon and Zerroukat, Mohamed and Thuburn, John},
  year = 2014,
  journal = {Quarterly Journal of the Royal Meteorological Society},
  volume = {140},
  number = {682},
  pages = {1505--1520},
  issn = {1477-870X},
  doi = {10.1002/qj.2235},
  copyright = {\copyright{} 2013 Royal Meteorological Society and Crown Copyright, the Met Office. Quarterly Journal of the Royal Meteorological Society \copyright{} 2013 Royal Meteorological Society},
  langid = {english}
}

@article{burgess_troposphere--stratosphere_2013,
    chapter = {Journal of the Atmospheric Sciences},
    title = {The {Troposphere}-to-{Stratosphere} {Transition} in {Kinetic} {Energy} {Spectra} and {Nonlinear} {Spectral} {Fluxes} as {Seen} in {ECMWF} {Analyses}},
    volume = {70},
    issn = {0022-4928, 1520-0469},
    url = {https://journals.ametsoc.org/view/journals/atsc/70/2/jas-d-12-0129.1.xml},
    doi = {10.1175/JAS-D-12-0129.1},
    language = {EN},
    number = {2},
    urldate = {2026-01-08},
    journal = {Journal of the Atmospheric Sciences},
    publisher = {American Meteorological Society},
    author = {Burgess, B. H. and Erler, Andre R. and Shepherd, Theodore G.},
    month = feb,
    year = {2013},
    pages = {669--687},
}

@article{waite_mesoscale_2009,
    chapter = {Journal of the Atmospheric Sciences},
    title = {The {Mesoscale} {Kinetic} {Energy} {Spectrum} of a {Baroclinic} {Life} {Cycle}},
    volume = {66},
    issn = {0022-4928, 1520-0469},
    url = {https://journals.ametsoc.org/view/journals/atsc/66/4/2008jas2829.1.xml},
    doi = {10.1175/2008JAS2829.1},

    language = {EN},
    number = {4},
    urldate = {2026-01-08},
    journal = {Journal of the Atmospheric Sciences},
    publisher = {American Meteorological Society},
    author = {Waite, Michael L. and Snyder, Chris},
    month = apr,
    year = {2009},
    pages = {883--901},
}

@article{skamarock_vertical_2019,
    chapter = {Monthly Weather Review},
    title = {Vertical {Resolution} {Requirements} in {Atmospheric} {Simulation}},
    volume = {147},
    issn = {1520-0493, 0027-0644},
    url = {https://journals.ametsoc.org/view/journals/mwre/147/7/mwr-d-19-0043.1.xml},
    doi = {10.1175/MWR-D-19-0043.1},
    language = {EN},
    number = {7},
    urldate = {2026-01-08},
    journal = {Monthly Weather Review},
    publisher = {American Meteorological Society},
    author = {Skamarock, William C. and Snyder, Chris and Klemp, Joseph B. and Park, Sang-Hun},
    month = jul,
    year = {2019},
    pages = {2641--2656},
}

@article{izzo_strong_2022,
    title = {Strong {Stability} {Preserving} {Runge}–{Kutta} and {Linear} {Multistep} {Methods}},
    volume = {48},
    issn = {1735-8515},
    url = {https://doi.org/10.1007/s41980-022-00731-x},
    doi = {10.1007/s41980-022-00731-x},
    language = {en},
    number = {6},
    urldate = {2025-12-15},
    journal = {Bulletin of the Iranian Mathematical Society},
    author = {Izzo, Giuseppe and Jackiewicz, Zdzislaw},
    month = dec,
    year = {2022},
    pages = {4029--4062},
}

@article{skamarockFullyCompressibleNonhydrostatic2021,
    chapter = {Monthly Weather Review},
    title = {A {Fully} {Compressible} {Nonhydrostatic} {Deep}-{Atmosphere} {Equations} {Solver} for {MPAS}},
    volume = {149},
    issn = {1520-0493, 0027-0644},
    url = {https://journals.ametsoc.org/view/journals/mwre/149/2/MWR-D-20-0286.1.xml},
    doi = {10.1175/MWR-D-20-0286.1},
    language = {EN},
    number = {2},
    urldate = {2026-03-16},
    journal = {Monthly Weather Review},
    publisher = {American Meteorological Society},
    author = {Skamarock, William C. and Ong, Hing and Klemp, Joseph B.},
    month = feb,
    year = {2021},
    pages = {571--583},
}

@article{benacchioSemiimplicitSemiLagrangianModelling2016,
    title = {Semi-implicit semi-{Lagrangian} modelling of the atmosphere: a {Met} {Office} perspective},
    volume = {7},
    issn = {2038-0909},
    shorttitle = {Semi-implicit semi-{Lagrangian} modelling of the atmosphere},
    url = {https://reference-global.com/article/10.1515/caim-2016-0020},
    doi = {10.1515/caim-2016-0020},
    language = {English},
    number = {3},
    urldate = {2026-03-16},
    journal = {Communications in Applied and Industrial Mathematics},
    publisher = {Italian Society for Applied and Industrial Mathemathics},
    author = {Benacchio, Tommaso and Wood, Nigel},
    month = oct,
    year = {2016},
    pages = {4--25},
}

@misc{FiredrakeUserManual,
  title        = {Firedrake User Manual},
  author       = {David A. Ham and Paul H. J. Kelly and Lawrence Mitchell and Colin J. Cotter and Robert C. Kirby and Koki Sagiyama and Nacime Bouziani and Sophia Vorderwuelbecke and Thomas J. Gregory and Jack Betteridge and Daniel R. Shapero and Reuben W. Nixon-Hill and Connor J. Ward and Patrick E. Farrell and Pablo D. Brubeck and India Marsden and Thomas H. Gibson and Miklós Homolya and Tianjiao Sun and Andrew T. T. McRae and Fabio Luporini and Alastair Gregory and Michael Lange and Simon W. Funke and Florian Rathgeber and Gheorghe-Teodor Bercea and Graham R. Markall},
  organization = {Imperial College London and University of Oxford and Baylor University and University of Washington},
  edition      = {First edition},
  year         = {2023},
  month        = {5},
  doi          = {10.25561/104839},
}

@article{kowalczykCompatibleFiniteElements2025,
    title = {On compatible finite elements for atmosphere modelling},
    url = {https://hdl.handle.net/10044/1/127314},
    doi = {10.25560/127314},
    language = {en},
    urldate = {2026-03-17},
    author = {Kowalczyk, Karina},
    month = jun,
    year = {2025},
}

@misc{chewOnestepBlendedSoundproofcompressible2021,
    title = {A one-step blended soundproof-compressible model with balanced data assimilation: theory and idealised tests},
    shorttitle = {A one-step blended soundproof-compressible model with balanced data assimilation},
    url = {https://arxiv.org/abs/2103.11861v2},
    doi = {10.1175/MWR-D-21-0175.1},
    language = {en},
    urldate = {2026-03-17},
    journal = {arXiv.org},
    author = {Chew, Ray and Benacchio, Tommaso and Hastermann, Gottfried and Klein, Rupert},
    month = mar,
    year = {2021},
    doi = {10.1175/MWR-D-21-0175.1},
}
\end{document}